%% file: main.tex
\documentclass[a4paper, twoside]{report}

\usepackage[english]{babel}
\usepackage[utf8]{inputenc}
\usepackage[T1]{fontenc}

\usepackage[a4paper,top=3cm,bottom=2cm,left=3cm,right=3cm,marginparwidth=1.75cm]{geometry}

\usepackage[sorting=none]{biblatex}
\usepackage{afterpage}
\usepackage{amsmath}
\usepackage{amsthm}
\usepackage{csquotes}
\usepackage{enumitem}
\usepackage{graphicx}
\usepackage{grffile}
\usepackage{lipsum}
\usepackage{booktabs}
\usepackage{float}

\usepackage{listings}
\usepackage[labelfont=bf,justification=centering]{caption}

\usepackage[ruled,linesnumbered,vlined]{algorithm2e}
\makeatletter
\renewcommand{\SetKwInOut}[2]{%
  \sbox\algocf@inoutbox{\KwSty{#2}\algocf@typo:}%
  \expandafter\ifx\csname InOutSizeDefined\endcsname\relax
    \newcommand\InOutSizeDefined{}\setlength{\inoutsize}{\wd\algocf@inoutbox}%
    \sbox\algocf@inoutbox{\parbox[t]{\inoutsize}{\KwSty{#2}\algocf@typo:\hfill}~}\setlength{\inoutindent}{\wd\algocf@inoutbox}%
  \else
    \ifdim\wd\algocf@inoutbox>\inoutsize%
    \setlength{\inoutsize}{\wd\algocf@inoutbox}%
    \sbox\algocf@inoutbox{\parbox[t]{\inoutsize}{\KwSty{#2}\algocf@typo:\hfill}~}\setlength{\inoutindent}{\wd\algocf@inoutbox}%
    \fi%
  \fi
  \algocf@newcommand{#1}[1]{%
    \ifthenelse{\boolean{algocf@inoutnumbered}}{\relax}{\everypar={\relax}}%
    {\let\\\algocf@newinout\hangindent=\inoutindent\hangafter=1\parbox[t]{\inoutsize}{\KwSty{#2}\algocf@typo:\hfill}~##1\par}%
    \algocf@linesnumbered
  }}%
\makeatother

\usepackage{bm}

\usepackage{fancyhdr}
\pdfoutput=1 
\usepackage{tipa}
\usepackage{tcolorbox}
\definecolor{mdgrey}{rgb}{0.8, 0.8, 0.8}
\usepackage[framemethod=tikz]{mdframed}
\usepackage{lipsum}
\newtheoremstyle{defi}
  {\topsep}%
  {\topsep}%
  {\normalfont}%
  {}%
  {\bfseries}%
  {:}%
  {.5em}%
  {\thmname{#1}\thmnote{~(#3)}}%
\theoremstyle{defi}
\newmdtheoremenv{definitioni}{Definition}
\newmdtheoremenv[
hidealllines=true,
leftline=true,
innertopmargin=0pt,
innerbottommargin=0pt,
linewidth=4pt,
linecolor=gray!40,
innerrightmargin=0pt,
]{definitionii}{Definition}
\newmdtheoremenv[
roundcorner=5pt,
innertopmargin=0pt,
innerbottommargin=5pt,
linewidth=4pt,
linecolor=gray!40,
]{definitioniii}{Definition}

\usepackage[colorinlistoftodos]{todonotes}
\usepackage[colorlinks=true, allcolors=blue]{hyperref}
\usepackage{subcaption}
\usepackage{algpseudocode}
\usepackage[linesnumbered,ruled]{algorithm2e}
\usepackage{amsmath}

\title{Maximum Return on Investment for a Domestic Photovoltaic Installation}

\author{Tom Nonnenmacher}

\bibliography{bibs/bibliography.bib}

\begin{document}
\input{title/title.tex}

\begin{abstract}
The rising energy prices in Europe and the urgent need to address global warming have sparked a significant increase in the installation of domestic photovoltaic systems to harness solar energy. However, since solar energy is available only during daytime hours and its availability varies daily, effectively shifting energy use becomes crucial. Whilst batteries can assist in storing excess energy, their high prices hinder their widespread adoption. In this study, we explore the importance of load to maximise return on investment.

We propose an incremental approach to fitting load profiles into the production envelope, allowing for practical implementation.
We compare different meter resolutions: 1 second, 5 minutes, 15 minutes, and 1 hour. 
Our analysis reveals that making real-time decisions (per second) leads to significant energy savings of 16\% compared to hourly decisions. 

Furthermore, we explore three types of device management strategies:
\begin{enumerate}
\item ON/OFF management independent of PV production,
\item ON/OFF management based on the current PV production,
\item ON/OFF management based on both current and forecasted PV production, utilising an optimal fit algorithm.
\end{enumerate}
Through our study, we demonstrate that our implementation of the third approach outperforms a standard management approach, resulting in more than 17\% cost savings. This study provides insights into the optimisation of load-shifting strategies in domestic photovoltaic installations, highlighting the importance of load control and the potential benefits in maximising the utilisation of solar energy while minimising energy costs and environmental impact.
\end{abstract}
%

\tableofcontents


\input{introduction/introduction.tex}
\input{chapters/methodology.tex}
\input{chapters/results.tex}
\input{conclusion/conclusion.tex}

\printbibliography

\end{document}

%% file: title/title.tex
\begin{titlepage}

\newcommand{\HRule}{\rule{\linewidth}{0.5mm}} 


\includegraphics[width=8cm]{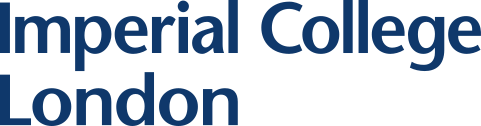}\\[1cm] 
 

\center 


\textsc{\LARGE Bachelor Individual Project}\\[1.5cm] 
\textsc{\Large Imperial College London}\\[0.5cm] 
\textsc{\large Department of Computing}\\[0.5cm] 

\makeatletter
\HRule \\[0.6cm]
{ \huge \bfseries \@title}\\[0.6cm] 
\HRule \\[1.5cm]
 

\begin{minipage}{0.4\textwidth}
\begin{flushleft} \large
\emph{Author:}\\
\@author 
\end{flushleft}
\end{minipage}
~
\begin{minipage}{0.4\textwidth}
\begin{flushright} \large
\emph{Supervisor:} \\
Prof. Jenny Nelson \\[1.2em] 
\emph{Co-supervisor:} \\
Dr. Benedict Winchester 
\end{flushright}
\end{minipage}\\[2cm]
\makeatother



{\large June 9, 2023}\\[2cm] 

\vfill 

\end{titlepage}

%% file: introduction/introduction.tex
\chapter{Introduction}

\section{Motivation}
Rising energy prices in Europe in 2022 and global warming causes the installation of 
domestic photovoltaic (PV) installations to harvest lost solar energy \cite{sgaravatti2022national, houghton2005global, kamat2007meeting}. Solar energy is available only during daytime hours and varies day by day. It is therefore key to shift energy use (load) to the time when solar power is available. Batteries can help to move energy availability but have a prohibitive high price at the time of the writing \cite{battery}. We take the perspective of domestic use, where the production and consumption of energy happen at home, thereby avoiding the cost and loss of transporting energy across distance.

Allocation of load of solar energy is paramount to optimise the use of the “free” solar energy.
We define typical loads of the home as:
\begin{itemize}
\item Electric car
\item Domestic hot water
\item Pool pump
\item Feed energy back into the grid (sell energy)
\item other (lighting, washing machines, cooking, PCs, etc.)
\end{itemize}

The problem becomes an optimisation problem that includes the learning of the parameters. We are trying to fit the various loads into the solar production curve of the day. At the time at which we make a decision for each and all loads, we do not have certainty about the coming solar energy available, and we do not know the energy use of each load. By using these assumptions we ensure that the practical implementation of our work fits universally and can be implemented without prior knowledge, or manual re-adjustment which itself is prone to error. 

\section{Practical considerations and context}

Our data from the 1st of January 2023 to the 10th of April 2023 is from a 14.7 kWp PV installation in the South of France. The time is in UTC. In our study, we will concentrate on developing algorithms tailored to this particular household, utilising its specific data and focusing on the control of three devices: the electric car, the pool pump, and the domestic hot water system.

France, as a highly industrialised nation, relies heavily on its energy infrastructure to meet the growing demands of its population and economy. The energy system in France predominantly comprises a mix of nuclear, renewable, and fossil fuel-based power generation sources \cite{mathew2022nuclear}. While nuclear energy has historically played a significant role in the country's energy production, there has been an increasing focus on incorporating renewable energy sources, such as solar, to diversify the energy mix and reduce carbon emissions \cite{lebrouhi2022energy}.

The varying availability of solar energy, coupled with fluctuating electricity demand patterns, necessitates the development of innovative strategies to optimise energy consumption and enhance the flexibility of the grid. Load-shifting, as a demand-side management technique, offers a promising solution to effectively balance electricity supply and demand while considering the constraints of the energy system.

\section{Literature overview}
Load shifting, a key aspect of demand-side management, has gained significant attention in recent years due to its potential to optimise energy consumption, reduce peak loads, and enhance overall grid efficiency. Load shifting involves adjusting the timing of energy usage, particularly in residential households, to take advantage of variable electricity prices, renewable energy availability, and grid demand patterns. Several studies have explored the application of load-shifting techniques in households equipped with energy storage units and electric mobility options, and have investigated the impact of demand-side management strategies on reducing household electricity consumption.

One notable research paper from Di Giorgio et al. presents a model predictive control framework for optimising load shifting in residential settings with energy storage systems \cite{di2012model}. The study focuses on developing strategies to effectively utilise the energy storage capacity to shift loads and minimise electricity costs.

Another research paper from Paetz et al. examines the load-shifting potentials in households that incorporate electric vehicles \cite{paetz2013load}. The study compares the behavior of users with the modeling results to evaluate the feasibility and effectiveness of load-shifting strategies in the context of electric mobility.

In the context of demand-side management, home appliance scheduling and peak load reduction play crucial roles in reducing household electricity consumption. Laicane et al. explore the impact of effective scheduling techniques on optimising energy usage and reducing overall electricity consumption in households \cite{laicane2015reducing}.

Furthermore, Fischer et al. focus on smart household usage profiling to recommend suitable energy tariffs and load-shifting strategies. The study utilises smart meter data and advanced profiling techniques to tailor load-shifting recommendations based on individual household energy consumption patterns \cite{fischer2013recommending}.

Lastly, Farzambehboudi et al. investigate the economic implications of load-shifting strategies in smart households \cite{farzambehboudi2018economic}. The study assesses the potential cost savings and benefits achieved through load-shifting techniques, shedding light on the economic feasibility and advantages of implementing such strategies.

In this paper, we aim to build upon the insights provided by these studies and contribute to the understanding of load-shifting methodologies in residential households. Our approach involves estimating future solar production, estimating each load profile, and optimally fitting load profiles within the production envelope.

%% file: chapters/methodology.tex
\chapter{Methodology}
In order to optimise the return on investment for our photovoltaic installation, we show our approach structured into three key sections: forecasting future solar production, estimating load profiles, and determining the optimal fitting of load profiles within the production envelope.

\section{The estimation of future solar production}
The PV output of each installation is different. Orientation, location, and obstacles vary the production profile. We want our results to be universally applicable anywhere in the world (southern or northern hemisphere, any location, any installation). Therefore, we proceed by measuring the solar production of the last 7 days and estimate future production. Each day, network connectivity issues and outages caused a few instances of missing data. To ensure data completeness, we utilised linear interpolation to fill in the gaps and obtained a total of 86400 data items per day, representing every second of the day. To estimate solar production, the starting point is the maximum possible production available for a given installation - the clear-sky model \cite{antonanzas2019clear}.  

First, the real production data from every day of January 2023 is being visualised. At first glance, the production curves resemble a normal distribution. We validate via the means squared error (MSE) how well the normal distribution approximates the maximum possible production (the clear sky model). If the normal distribution approximates well the clear sky model it will approximate also well the real production. We, therefore, start by approximating the real production data with a normal distribution to either validate or abandon the normal distribution as an approximation of the clear sky model.

We tried first to find the closest normal distribution for PV production of each given day. The normal distribution is defined by three parameters: the maximum, the mean, and the standard deviation. We selected the closest normal distribution for every production day in January 2023. 

\begin{figure}[H]
    \centering
    \includegraphics[width=9cm]{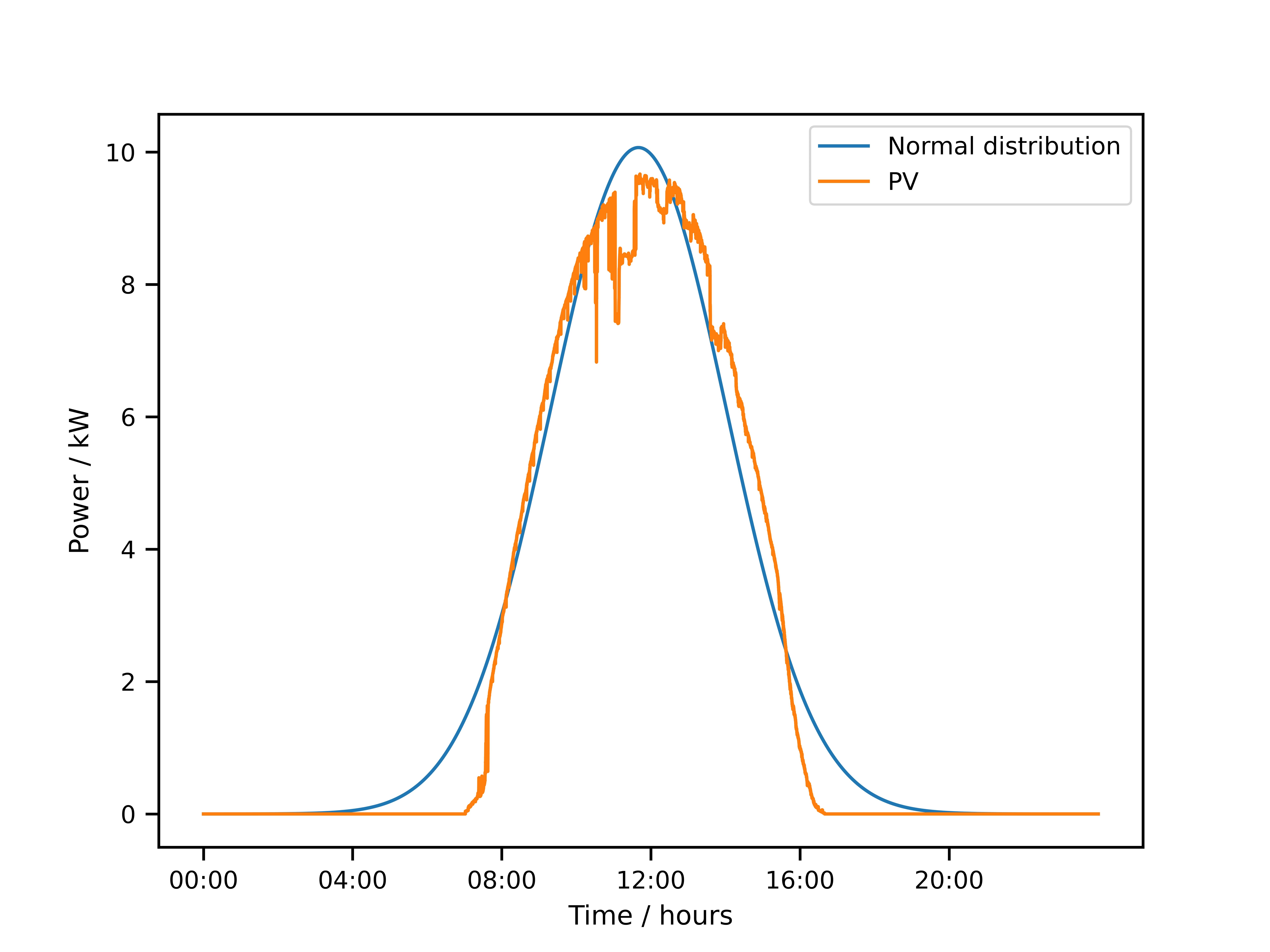}
    \caption{Normal distribution estimation for the PV production on the 31st of January 2023}
    \label{fig:normdist}
\end{figure}

To increase the speed of finding the closest distribution we optimised the search in the following way:
first, we selected the maximum value of the normal distribution in this range (in Watt):
\begin{equation} [max(PV) - 2000, max(PV) + 2000], \end{equation} second, we selected the mean of the normal distribution to be in the range from 8:20 am to 4:40 pm, and third the standard deviation in the range of from 4000 seconds to 20000 seconds. Within this search space, we chose the normal distribution with the smallest Mean Squared Error (MSE). 
Even the best normal distribution didn’t well approximate the production data, as we see on the following graph. The Figure \ref{fig:normdist} shows the day with the smallest MSE between a normal distribution and production data (31st January 2023). The reader sees immediately that the start and the end of the PV production are far overestimated. The MSE, computed using the W unit, is 362893.

We conclude that the normal distribution is not a good approximation of the production data, and therefore is also not a good approximation for the clear sky model. 

In the paper \cite{peratikou2022estimating} a method of estimating clear-sky PV electricity production without exogenous data is given. This article adopts a direct forecasting approach and thus, the proposed data-driven approach is dependent solely on historical time-series power output data from grid-connected without having to rely on the technical characteristics of the PV or on other energy or meteorological-related data, equipment, models, sky/satellite images, etc. 

\begin{figure}[H]
    \centering
    \includegraphics[width=9cm]{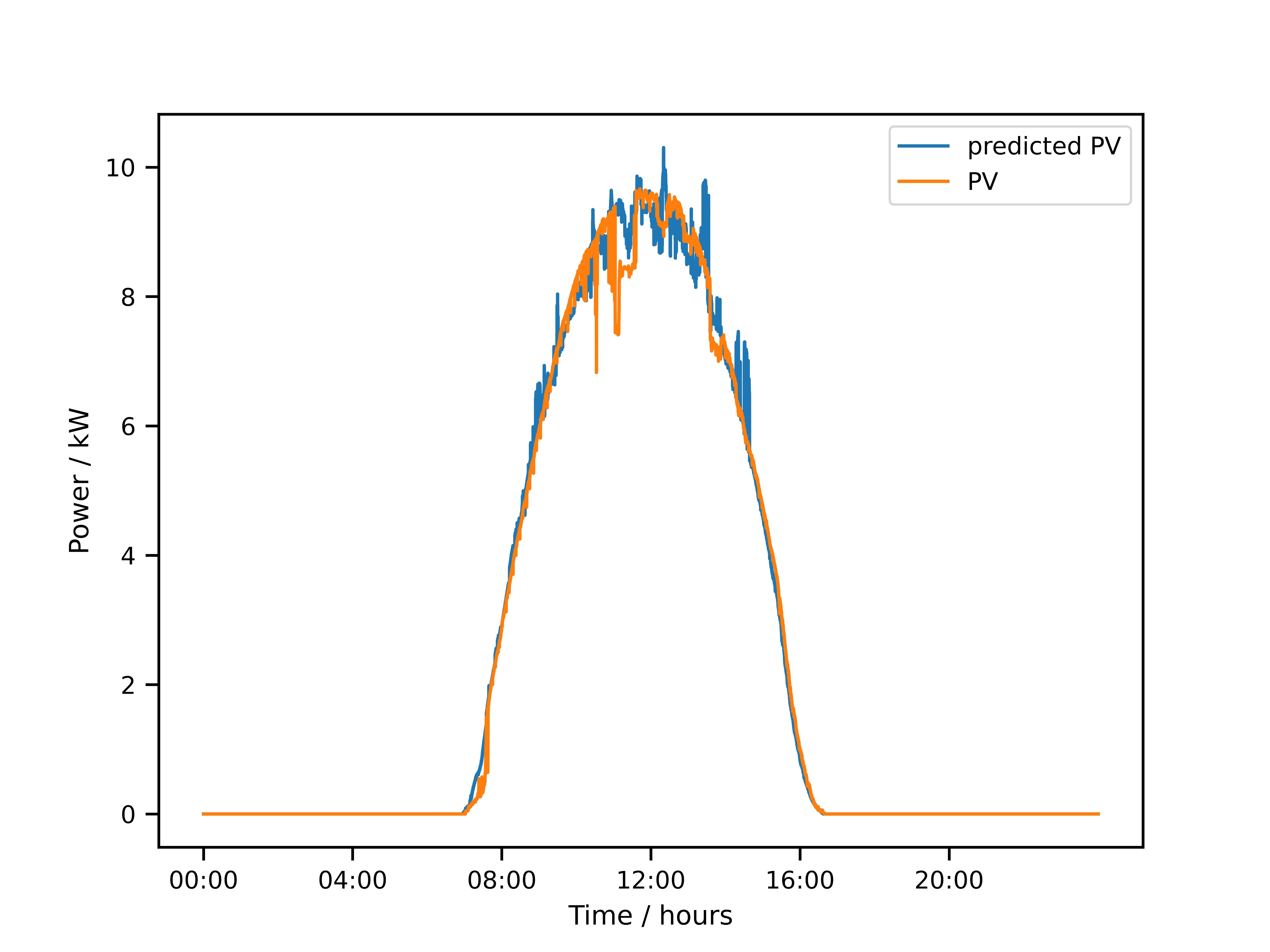}
    \caption{PV production forecasted and PV production on the 31st of January 2023}
    \label{fig:forecast}
\end{figure}

Their method to estimate the clear sky is following:
\begin{enumerate}
\item The quantities that characterise the shape and the distribution of PV power output time series at a given time are calculated first. The three criteria used are the mean value, the maximum value, and the standard deviation of the PV power output time series.
\item The calculated quantities were used as limits to the PV power output at each second. All the PV power output values outside the limits were excluded. Different limits were studied and tested in their study.
\item The predicted clear-sky signal is calculated by the resulting average of the PV power output time series within the limits for the corresponding month.
\end{enumerate}
The best method (with the smallest MSE) describes: 
\begin{equation}
\begin{split}
upper limit &= A + 1.5 SD \\
lower limit &= B. 
\end{split}
\end{equation}
A is the mean value of the examined PV power output time series at a specific time. \\
SD is the Standard Deviation of the examined PV power output time series at a specific time. \\
B is the average value of PV power output, excluding all points outside the upper limit and the lower limit band. \\
If we use the previous example of the 31st of January, with the use of the PV production of the last seven days (from the 24th of January 2023 to the 30th of January 2023), the algorithm outputs a predicted production for the 31st of January 2023 with a sufficient MSE of 68340, which is an 81\% decrease compared to the MSE of the normal distribution method (that knows the current PV production). The result is shown in Figure \ref{fig:forecast}.

\begin{figure}[H]
\subfloat{\includegraphics[width = 4cm]{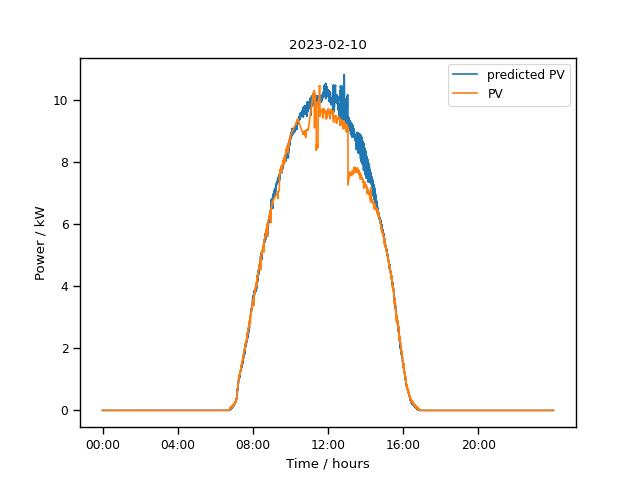}} 
\subfloat{\includegraphics[width = 4cm]{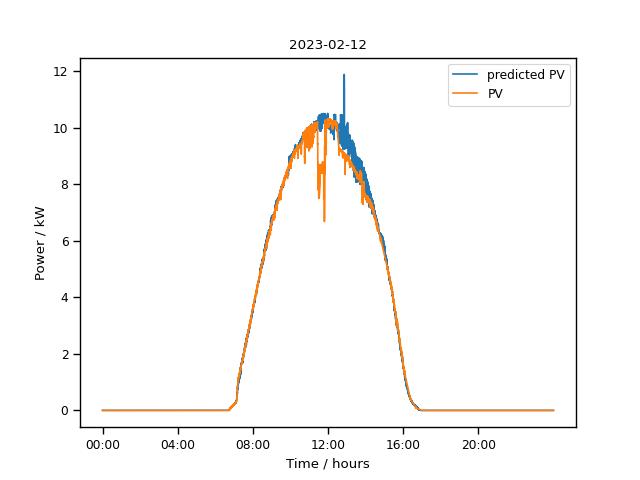}} 
\caption{The 10 best PV production forecasts from the 8th January to the 10th April 2023}
\label{fig:10best}
\end{figure}

To measure the performance of the PV forecasting, we measured all the clear sky data we had: from the 8th of January (because we need 7 past days) to the 10th of April 2023. 
We computed the mean per second of the difference between the real PV production and the PV production forecasted for every day. The result is:
the forecast of the PV production overestimates the real PV production per second with a mean of \textbf{2.12\%}.
Please find the following graphs for the 10 best PV forecasts
with \ref{fig:10best}. 

\begin{figure}[H]
\subfloat{\includegraphics[width = 4cm]{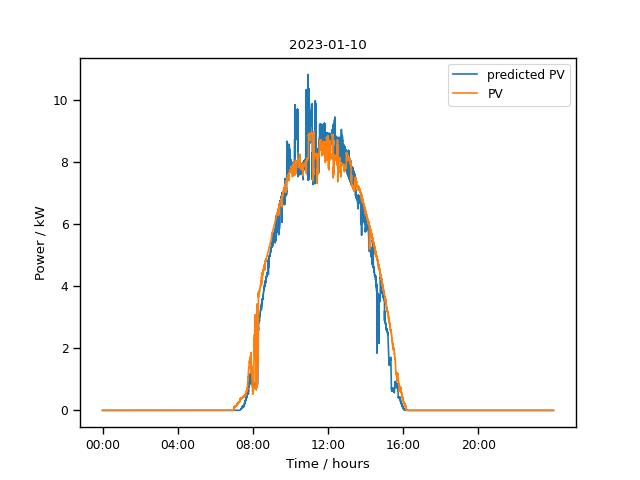}}
\subfloat{\includegraphics[width = 4cm]{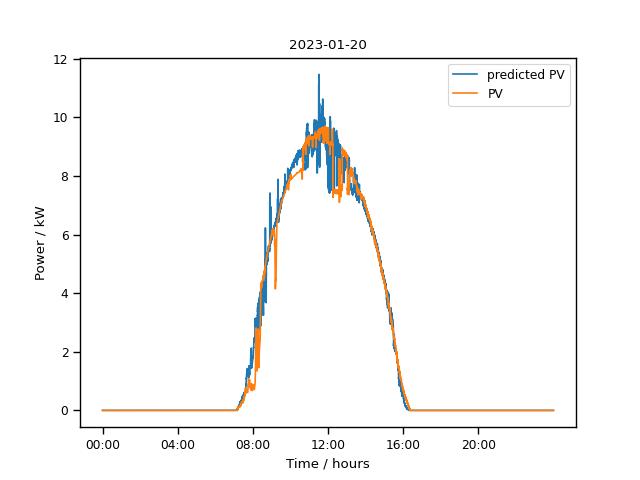}} 
\subfloat{\includegraphics[width = 4cm]{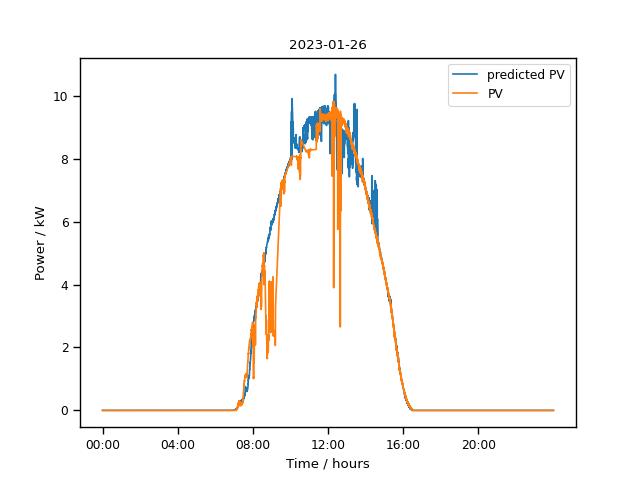}} 
\subfloat{\includegraphics[width = 4cm]{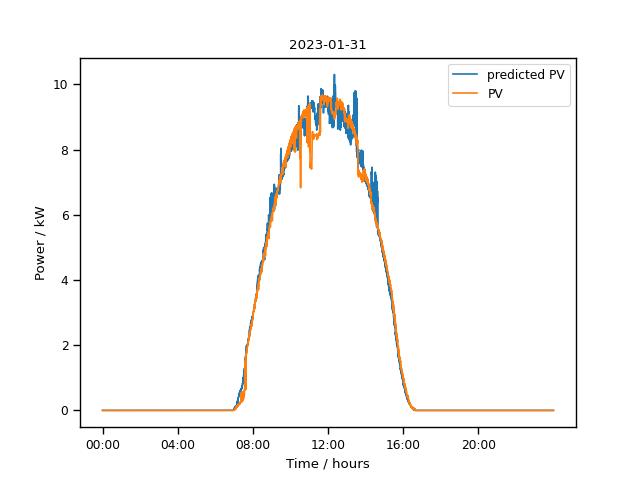}}
\end{figure}
\begin{figure}[H] 
\subfloat{\includegraphics[width = 4cm]{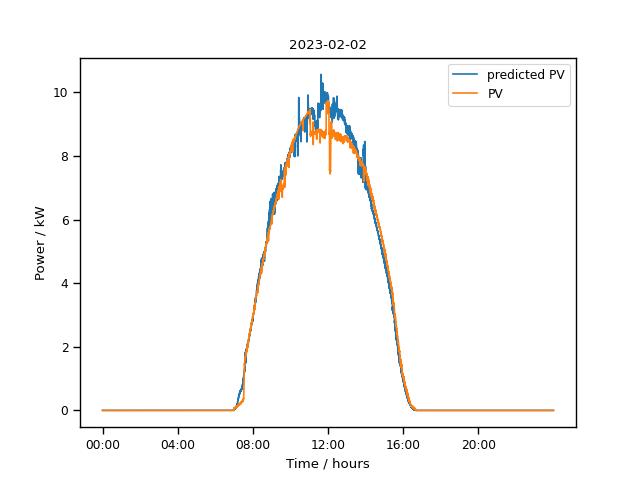}}
\subfloat{\includegraphics[width = 4cm]{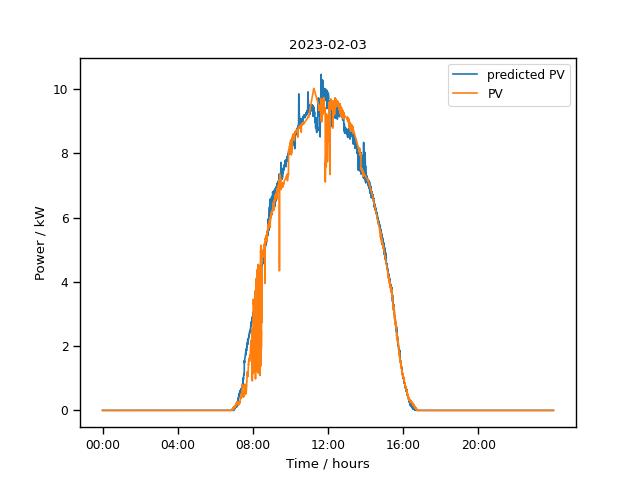}} 
\subfloat{\includegraphics[width = 4cm]{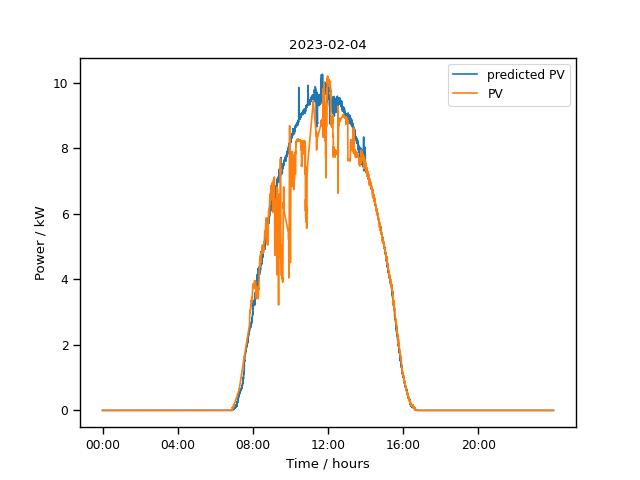}} 
\subfloat{\includegraphics[width = 4cm]{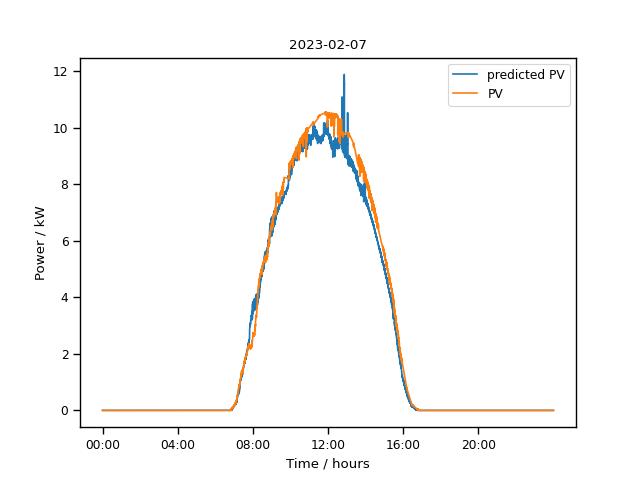}} 
\end{figure}
\clearpage

\section{The estimation of each load profile}
We assume that we have control over the power we give to each device and we can observe only the total energy load of all consumer loads in the home. We filter and estimate the load of a device by observing the total consumption of energy and the impact when we switch the device ON and OFF. By doing this multiple times, during multiple days we are able to distill a load profile within reasonable confidence intervals.
\subsection{The different controllable devices}

\subsubsection{Electric car} 
The majority of the electric-car charging is at constant power, e.g. cars come with a charger that allows for household plug charging (2.3 kW - 10 A at 230 V), or with a dedicated Wallbox at a given power, oftentimes this is 7.36 kW or 11 kW (corresponding to 32A, 48 A, respectively with 230 V) \cite{batterycapacity}. However, charging at very high power (fast charging) reduces the battery longevity \cite{al2022fast}.

In France, the average consumption is 50 minutes per person and per day, for an average distance of 36km \cite{art:1}. An electric car uses in average 0.161 kWh per km \cite{MGEV}. Therefore in this study, we assume the daily power requirement for the car is: 5.796 kWh. 

Hence, given this information, we fix the charging time for the car at \textbf{2 hours 30 minutes}, because it is the charging time needed to satisfy the average consumption in France if we use a household plug charging of 2.3 kW. Practically, a longer trip with the car can be foreseen and may require manual intervention to boost the electric-car battery to maximum capacity. 

\subsubsection{Domestic hot water} 
There are a vast amount of water boilers deployed in homes today, and many use electricity as the primary source to heat hot water until the target temperature is reached (typically 60 °C). When the target temperature is reached the boiler's security system shuts the heating OFF, and no electricity is consumed anymore beyond that point.

Per person and year, 800 kWh is the required energy for hot water in a household \cite{art:3}.
In France, there are an average of 2.19 people per household \cite{art:4}. So the estimated daily energy used in the household for hot water is 4.8 kWh.

Therefore, an average 2 kW water boiler requires 2.4 h of heating daily with 2 kWh of electricity.
During the day, this form of water heating allows for load shifting as showers and hot water consumption is concentrated at the beginning and the end of the day.

Usually, the typical pattern is full charge for 2 and half hours followed by brief intervals of full charging for a few minutes to maintain the water temperature at the desired level. In this study, we only include the main big initial heating period of \textbf{2 hours 30 minutes} and neglect the small periods.

\subsubsection{Pool pump} 
A pool pump is a device used to circulate and filter water in a swimming pool. It is responsible for pulling water from the pool, passing it through a filtration system to remove debris and contaminants, and then returning the clean water back to the pool.
The rule of thumb used for the running time of a pool pump that runs water through the filtration system is water temperature divided by 2, the result is converted in hours. 

So at 20 degrees water temperature, the pool pump is supposed to run for \textbf{10 hours} every day. Pool filtration and solar power go nicely hand-in-hand.

\subsection{Estimation of each device}
\label{sec:algo_load}
After having computed the PV forecasting, we need to estimate the consumption of a device. For perfect device management, the highest level of automation, and the avoidance of human errors, we use the automatically estimated device profiles to do the device and load management in the next chapter.

We want to estimate the consumption of a device, we have as input the time when the device is switched ON and the time when the device is switched OFF. Additionally, we have access to the global load of all devices during these periods. We assume the most common case, which is the presence of only one energy meter for the household. The problem can therefore be stated as: 
\begin{center}
    \textbf{Separating the energy consumed by one device from the sum of energy consumed by all devices.}
\end{center}

\begin{figure}[H]
\hfill
\subfloat{\includegraphics[width=7cm]{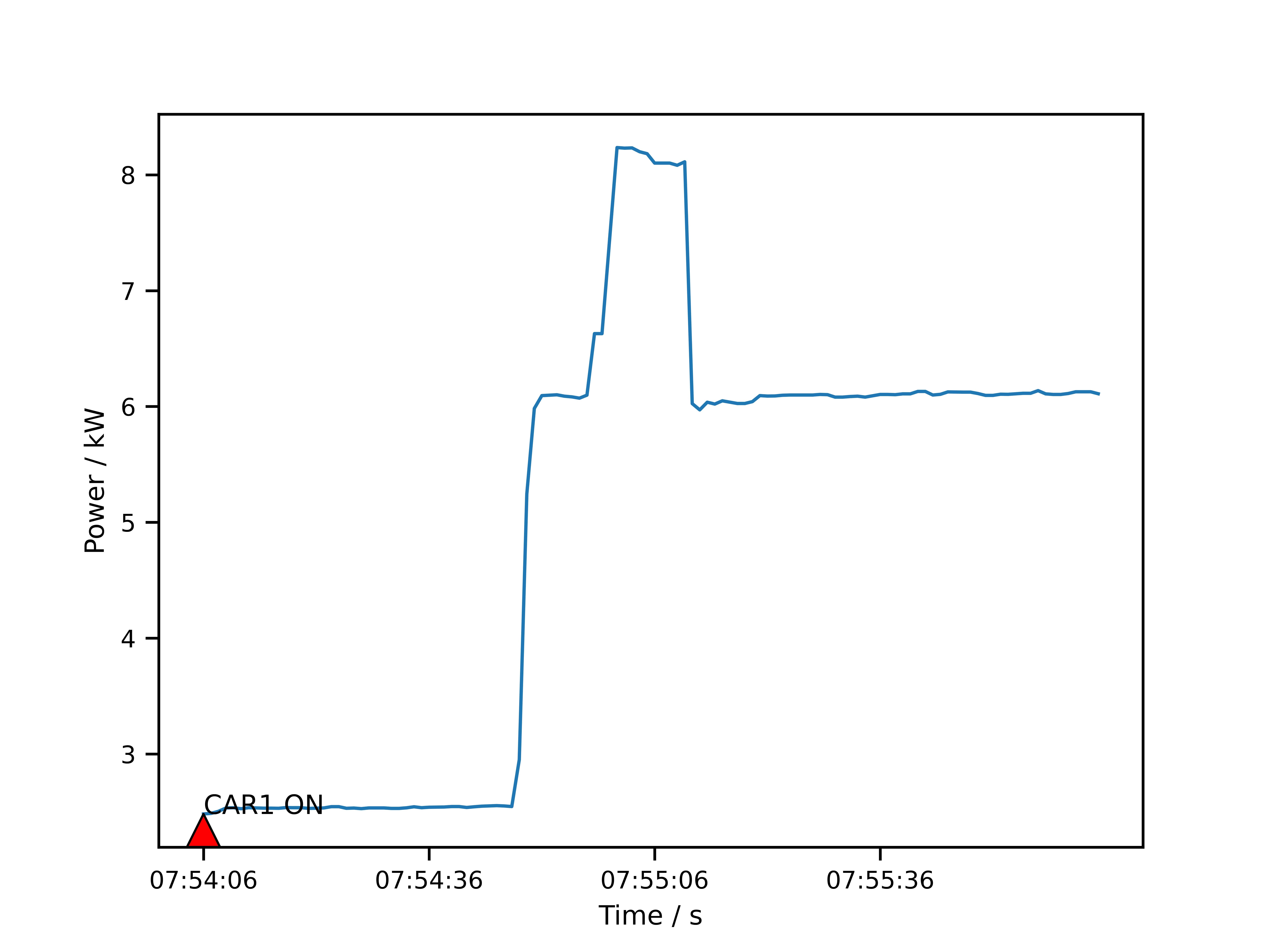}}
\hfill
\subfloat{\includegraphics[width=7cm]{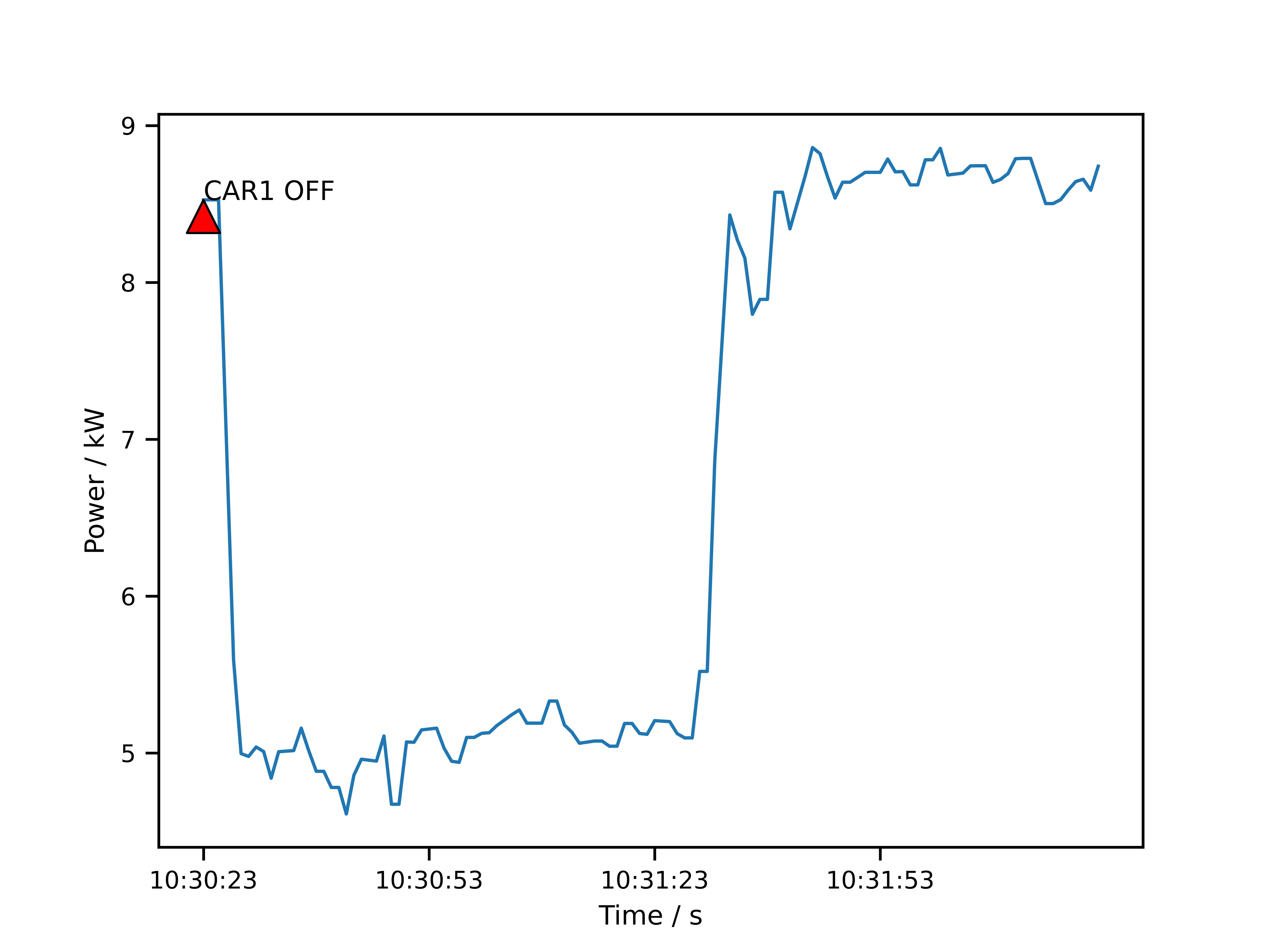}}
\hfill
\caption{The first ON and OFF of the electric CAR on the 1st April of 2023}
\label{fig:firstON}
\end{figure}

With the available ON/OFF timing information and the meter reading at all times, we search for the best estimate of the device consumption. We separate 2 cases:
\begin{itemize}
\item when the device switches ON, and
\item when the device switches OFF.
\end{itemize}
When the device switches ON, the load increases with a magnitude corresponding approximately to the consumption of the device. However, sometimes, the change is not instant, for example with the car (\ref{fig:firstON}). It often takes up to one minute after switching ON until the load becomes active and apparent on the meter reading.

When the device switches OFF, the load decreases with a magnitude corresponding approximately to the consumption of the device. The response is quicker after an OFF event than an ON event (see Figure \ref{fig:firstON}).

To find the consumption of the device we do edge detection on the global load with a window of 10 seconds: 
we shift a window of 10 seconds over the load, and we save the absolute value of the difference between the maximum and minimum inside the window as the device loads. To determine if it is an ON or an OFF event we verify the order, i.e. the minimum needs to be after the maximum for the ON event, while for the OFF event, it is the inverse. 

For an ON event, we search for the maximum load difference within a two-minute interval starting from the time when the device is switched ON (t to t + two minutes). Similarly, for an OFF event, we search for the maximum load difference within a one-minute interval starting from the time when the device is switched OFF (t to t + one minute). Furthermore, we apply exponential weighting to the window, meaning that an absolute difference measured between t and t + size of the window carries twice the weight of the absolute difference measured between t + size of the window and t + two times the size of the window. In summary, the consumption of a device is determined by the maximum weighted absolute difference, following the aforementioned conditions.

\subsection{Turbulent periods and base load}
\label{sec:calm_periods}
A "calm" period is a 10 minutes period of load that has a standard deviation corresponding to the 10\% of the lowest standard deviations of all the 10 minutes periods of load of any given day. A turbulent period is defined as a period that is not a "calm" period.
As we know, devices are sometimes switched ON/OFF during turbulent periods and thus can lead to wrong results. 
In fact, two devices can be launched at the same time, and it will be very difficult to find the consumption of each device.  To enhance the accuracy of our device estimation, we considered the possibility of measuring the load change when activating a device during regular "calm" periods, aiming for the highest level of precision possible.
For $n$ devices, we choose \begin{math} n\times(2+1)+1 \end{math} minutes to define a calm period. Three minutes per device (composed of two minutes for the ON period, and one minute for the OFF period) and one minute as a safety buffer for all devices. 

\begin{figure}[H]
\subfloat{\includegraphics[width = 5cm]{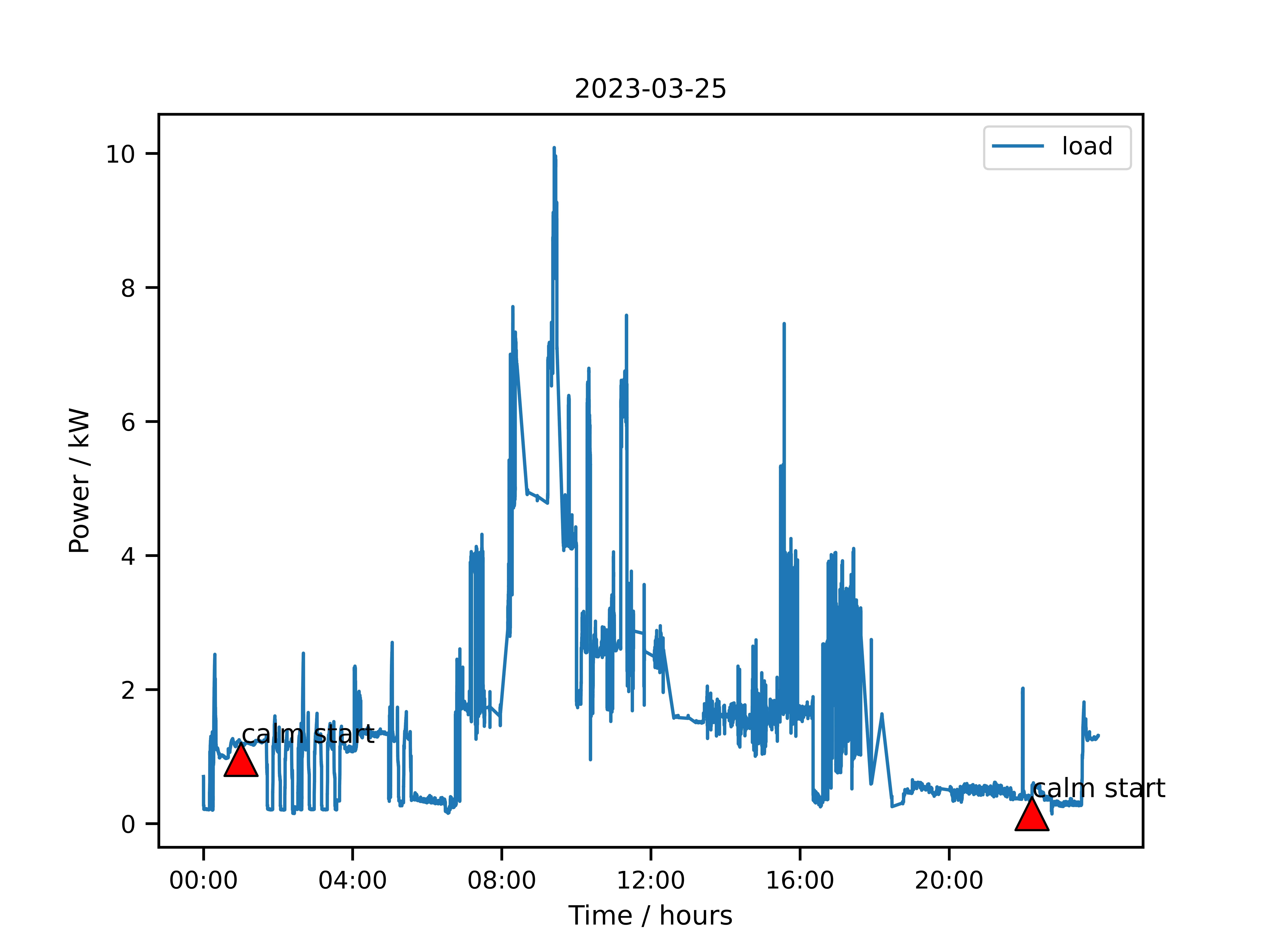}} 
\caption{The calm periods found from the 19th March 2023 to the 25th March 2023}
\label{fig:pastDays}
\end{figure}

To forecast a calm period in a day, we compute all of the calm periods during the past 7 days and require that at least 5 of them overlap. 
All these points forecasted the start of the calm periods for our current day. 
The graphs in Figure \ref{fig:pastDays} show that the two calm periods forecasted for the 26th of March 2023 start at 01:00 am and 10:22 pm, based on the data of the last 7 days (from the 19th of March 2023 to the 25th of March 2023). Therefore, on the 26th of March 2023, the measurement and device load estimation happens from 01:00 am to 01:10 am and from 10:22 pm to 10:32 pm.
\begin{figure}[H]
\subfloat{\includegraphics[width = 5cm]{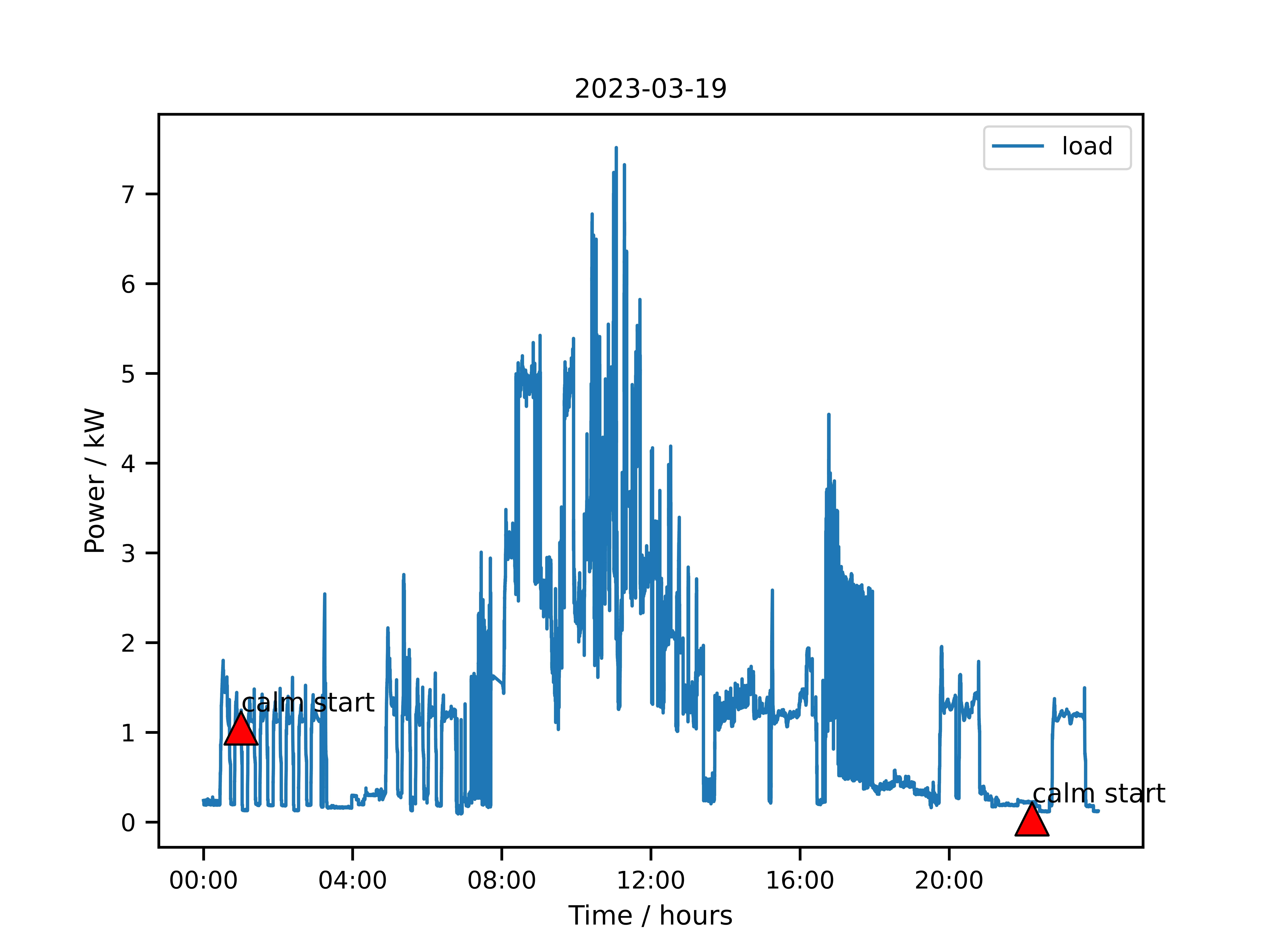}}
\subfloat{\includegraphics[width = 5cm]{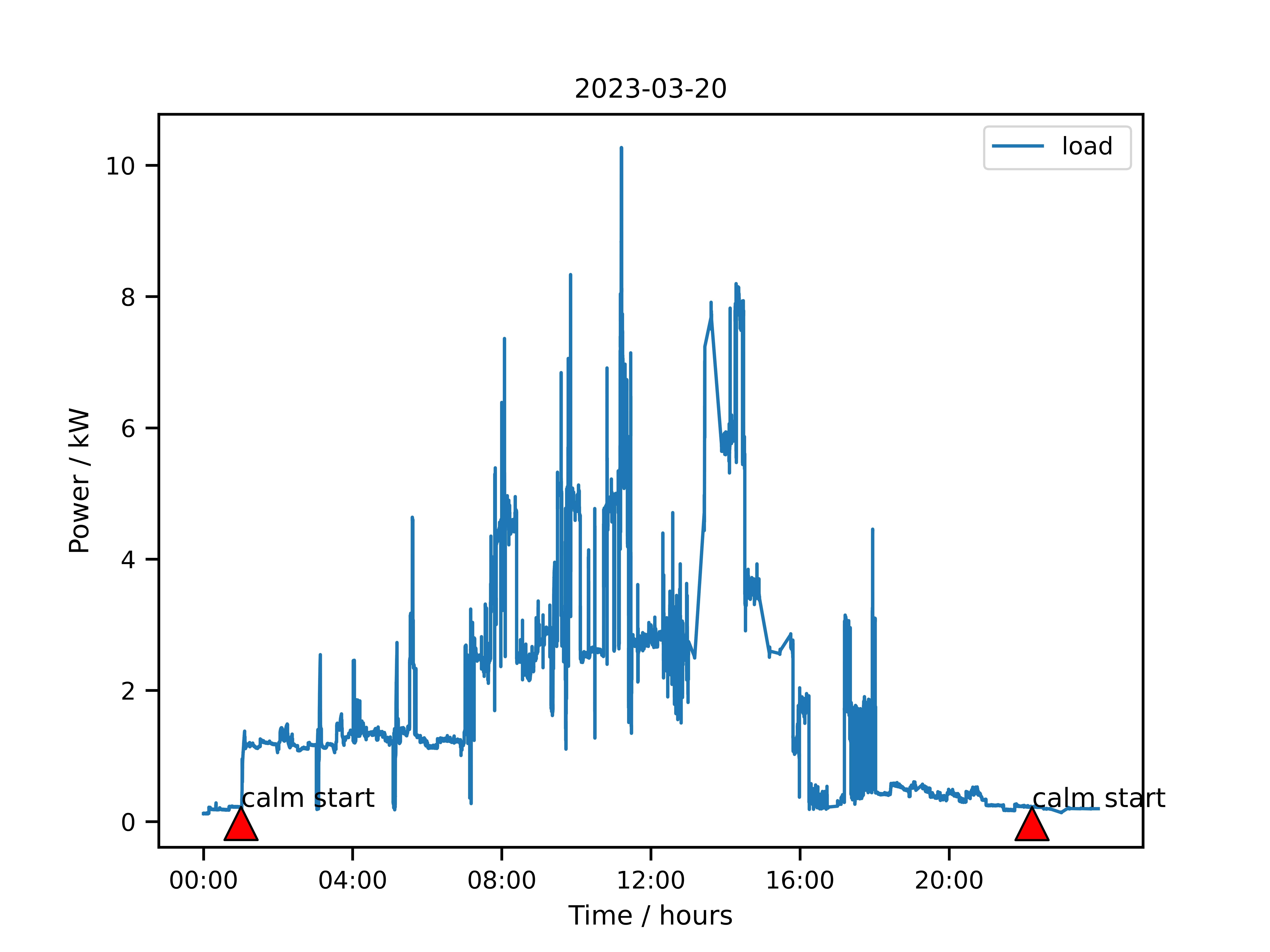}} 
\subfloat{\includegraphics[width = 5cm]{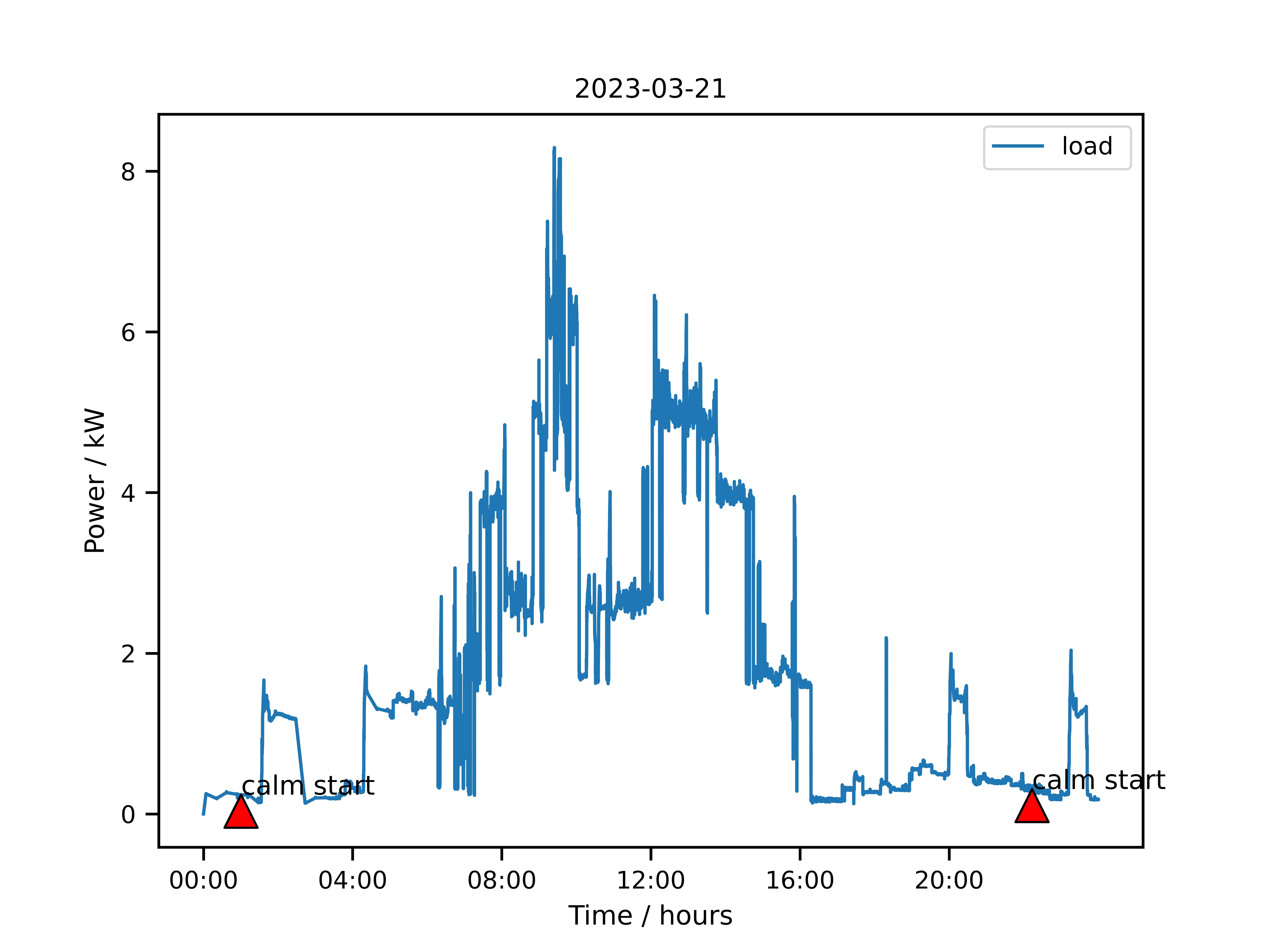}} 
\end{figure}
\begin{figure}[H]
\subfloat{\includegraphics[width = 5cm]{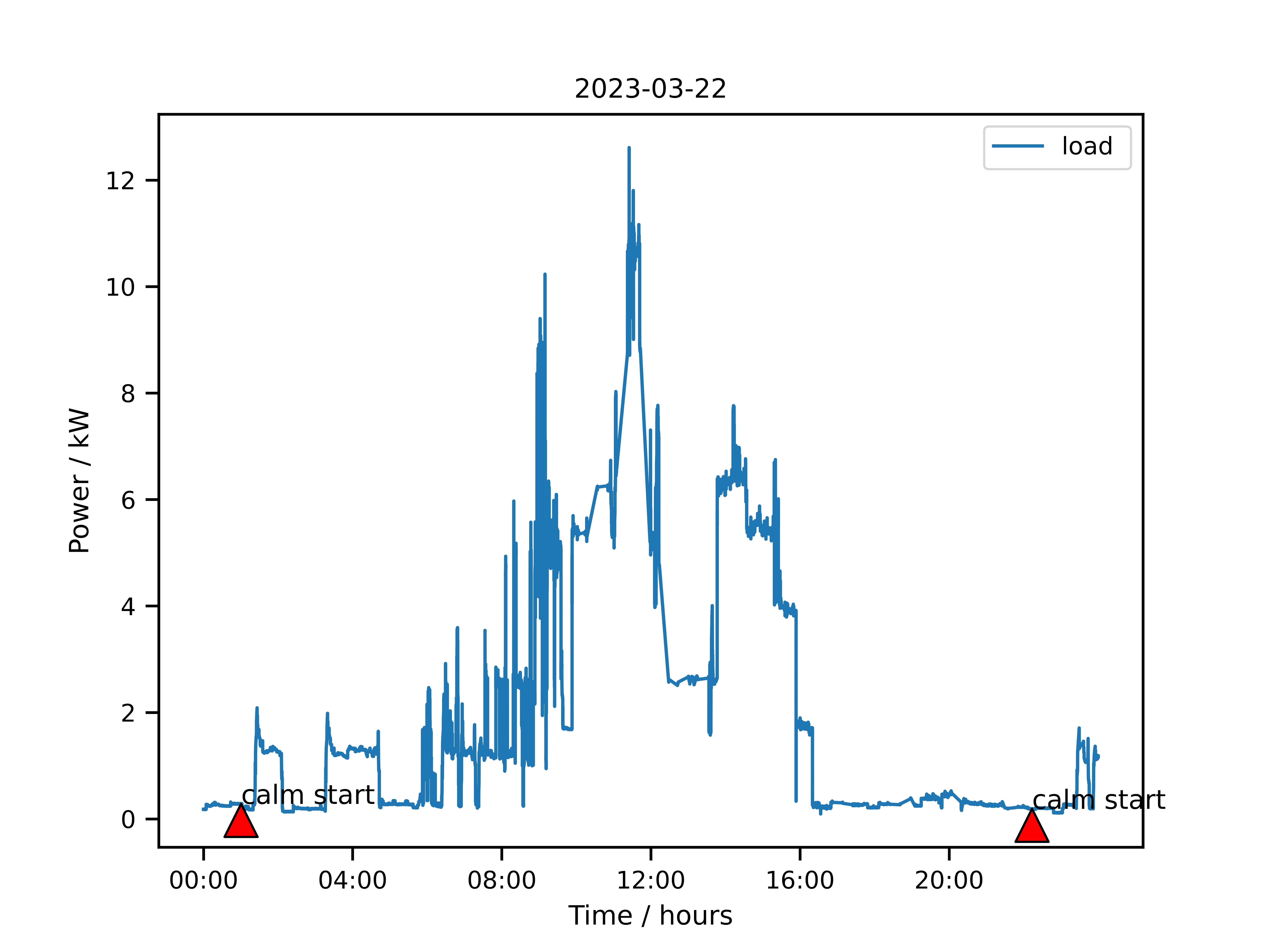}} 
\subfloat{\includegraphics[width = 5cm]{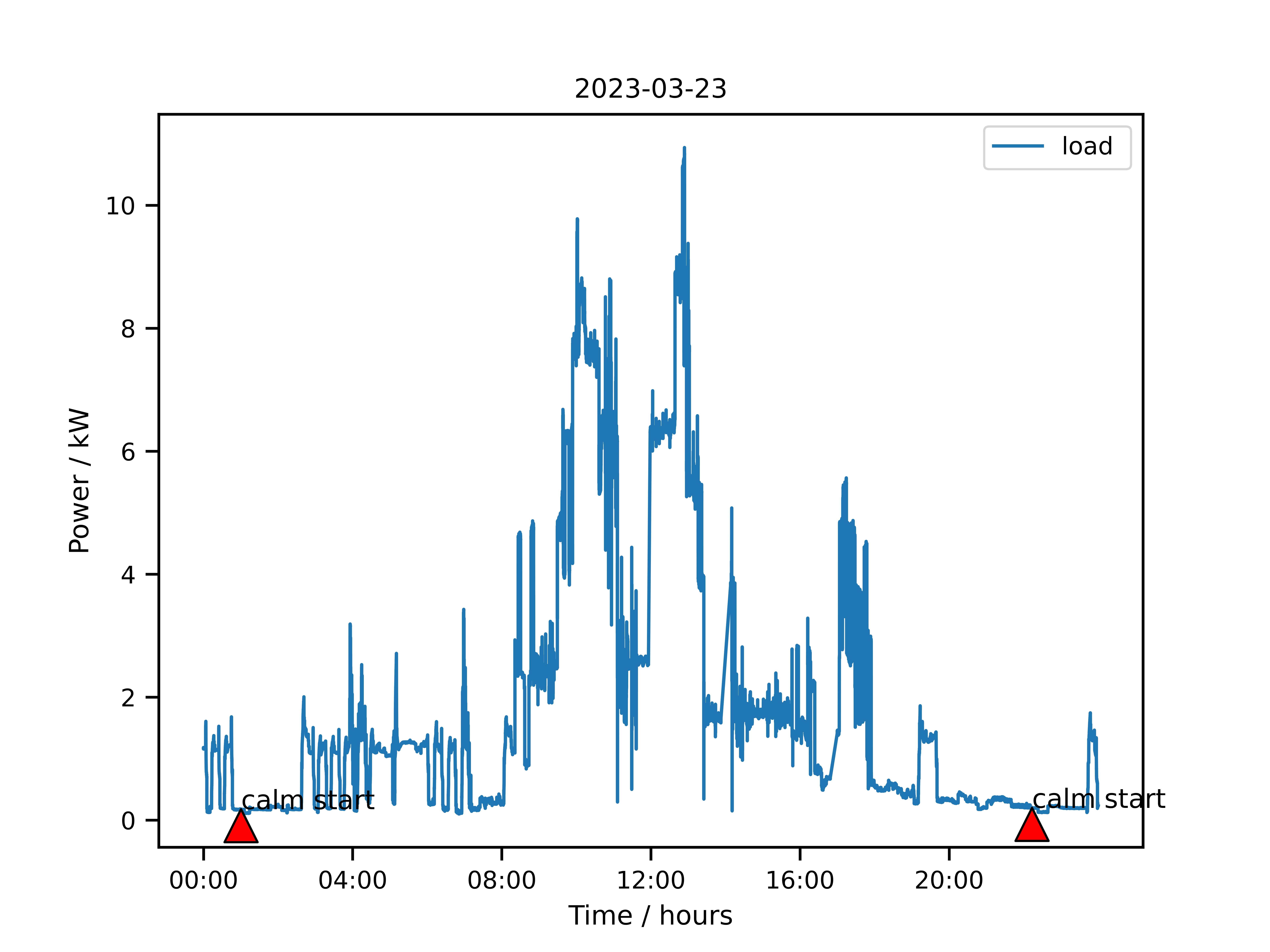}}
\subfloat{\includegraphics[width = 5cm]{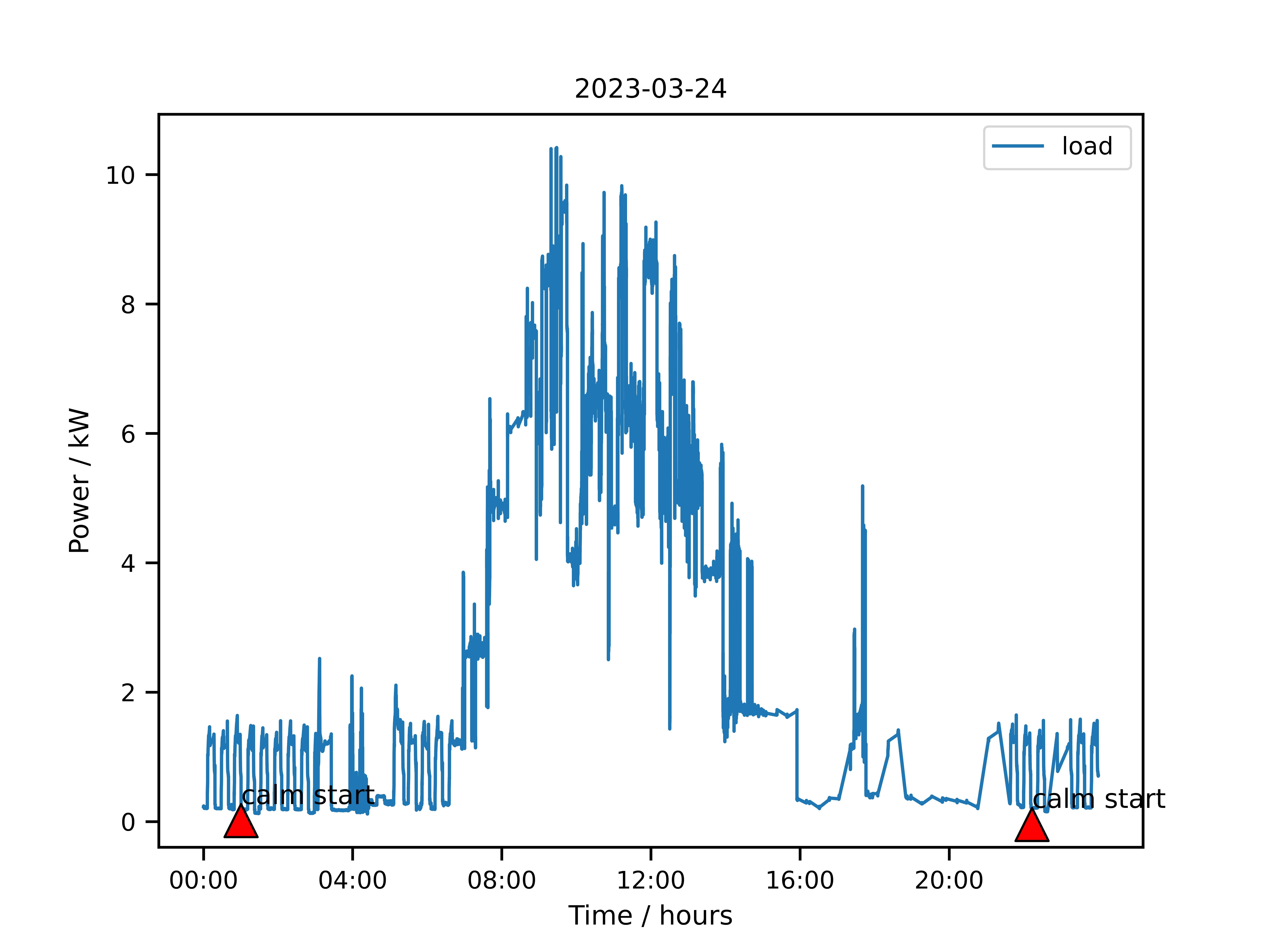}} 
\end{figure}
Calm periods are often during times outside of the PV production window (most of the calm periods found are during the night), which means that we will need to pay electricity from the grid to test and to estimate the device loads. Devices might change, therefore device estimation is done on a continuous basis, always using the information from the past seven days. Once seven days of data have been, we have enough data to even avoid testing outside the PV production window. By taking the median on the data measured, we avoid the outliers caused during turbulent periods and maintain good results. The algorithm with the seven past days found the correct device consumption of each device for every day by taking the median of all the consumption measured for ON/OFF events without using the "calm" periods.

\label{sec:base_load}
Then we needed to find the “base load” to simulate the load of the devices that we do not have control over, like the washing machine, the fridge, the lights, etc. The “base load” is the load when all the ON/OFF controllable devices are OFF. 
To compute it, we just computed the median of all the loads measured when all the controllable devices (the car, the pool pump, and the domestic hot water) were OFF. 

\section{Optimal fitting of load profiles into the production envelope}
We use an incremental approach and fit one additional load into an existing load of various other profiles into the production envelope. As the incremental approach can be repeated several times we have a practical method for all loads. The incremental deployment also reflects the reality of implementation. 

The problem now is to fit the load profile into the production envelope to maximise return on investment. First, money is saved by the consumption of cheap energy. Second, money is gained by selling excess capacity to the electricity provider.

Currently, the price to buy 1 kWh in France is 0.2€, and the price to sell 1 kWh to the electricity provider is 0.13€ \cite{energyprice}. Furthermore, we can not sell more than 6 kW per second.

With this information, we deduct the following constants for a day $j$:
\begin{align*} 
    E_{max,sell} &= 6 kWs  \\
    C_{e,sell} &= 0.13  EUR/kWh  \\
    C_{e,buy} &= 0.2  EUR/kWh 
\end{align*}
We name the 2 following functions for the load and the PV production: 
\begin{itemize}
\item P(j, t): "The PV production in kW for day $j$ at time $t$ seconds" 
\item L(j, t): "The load in kW for day $j$ at time $t$ seconds" \\
\end{itemize}
We, therefore, have one combined optimisation function of money that comes from saving and selling:
\\
\[
M(j) = S(j) - B(j)
\]
\begin{center}
which is the money in € for day $j$.
\end{center}

\[
S(j) =  \sum_{t=0}^{86399} min(E_{max,sell}, max(P(j, t) + L(j, t), 0))\times  \frac{C_{e,sell}}{3600s} 
\]
\begin{center}
which is the money in € we gain by selling energy on day $j$. 
\end{center}

\[
B(j) = \sum_{t=0}^{86399} max(0, -P(j, t) - L(j, t))\times  \frac{C_{e,buy}}{3600s} 
\]
\begin{center}
which is the money in € spent by buying energy on day $j$.
\end{center}

Our first idea was to work with basic conditions on thresholds (adaptative mode described in section \ref{sec:1}), for example, if the PV production deducted by the existing loads is greater than the consumption device \begin{math}i\end{math} at time \begin{math}t\end{math}, then we can switch device \begin{math}i\end{math} ON. If the PV production deducted by the existing loads is less than the consumption device \begin{math}i\end{math} at time \begin{math}t\end{math}, then we switch device \begin{math}i\end{math} OFF.

The problem with this algorithm is that it does not maximise electricity sales on sunny days, as sales have an upper bound. In fact, because it will turn on all the devices linearly with the beginning of the PV production, and then at the peak of the PV production, we can only sell 6 kW per second, so we will lose some PV production.

\begin{figure}[H]
    \centering
    \includegraphics[width=14.5cm]{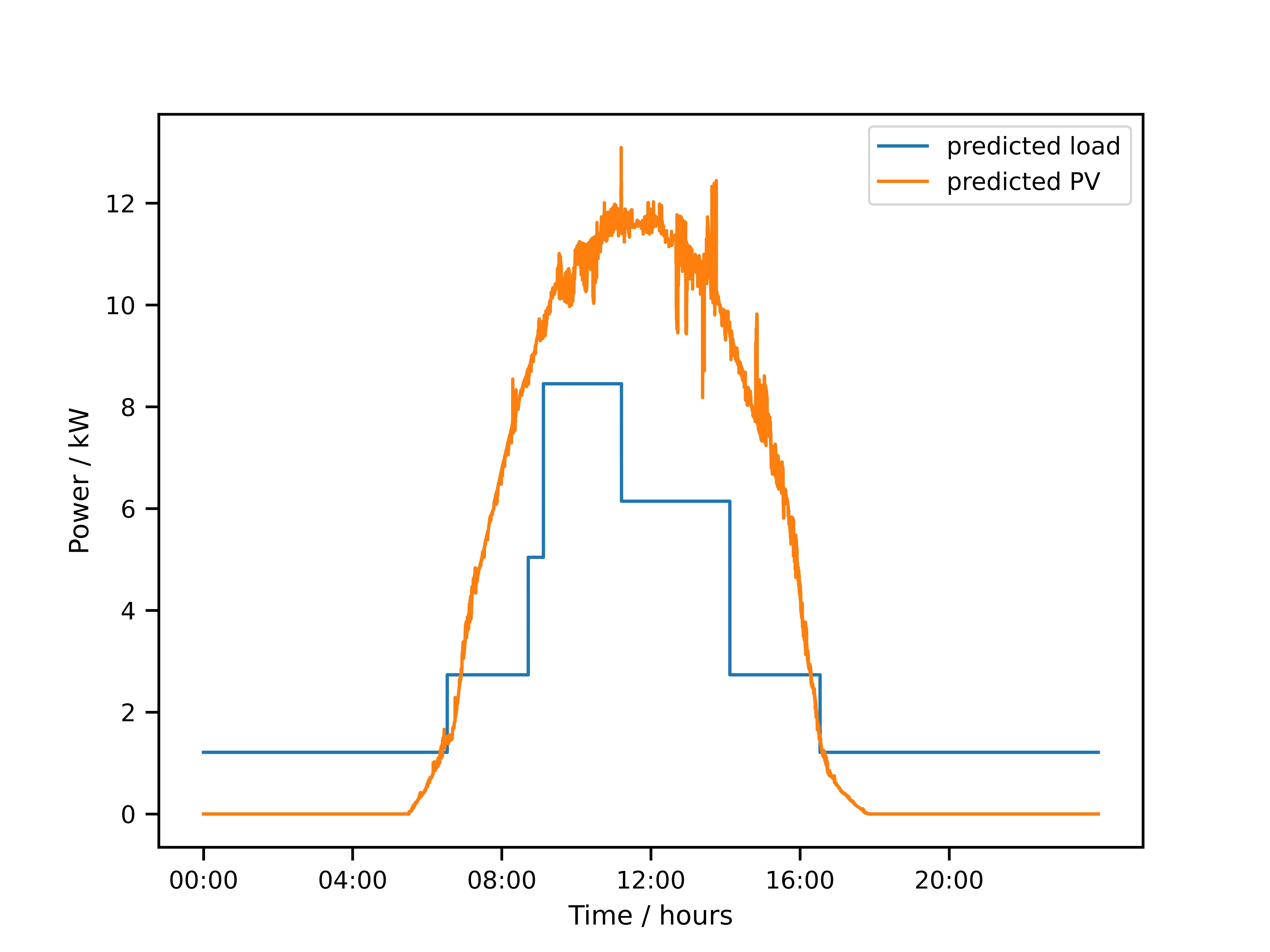}
    \caption{Forecasted load and PV production for the 1st April of 2023}
    \label{fig:rectangles}
\end{figure}

Then, we tried a different approach given the PV forecasting by maximising the function \begin{math}M\end{math}. The idea is to define the loads of each device as a rectangle. In fact, the height of the rectangle is the consumption of the device, and its length is the time the device needs to be on. We then place the rectangle of each device by the reverse order of length (rectangle with maximum length is placed first) to maximise the function \begin{math}M\end{math}, shown in Figure \ref{fig:rectangles}.

\begin{figure}[H]
    \centering
    \includegraphics[width=14.5cm]{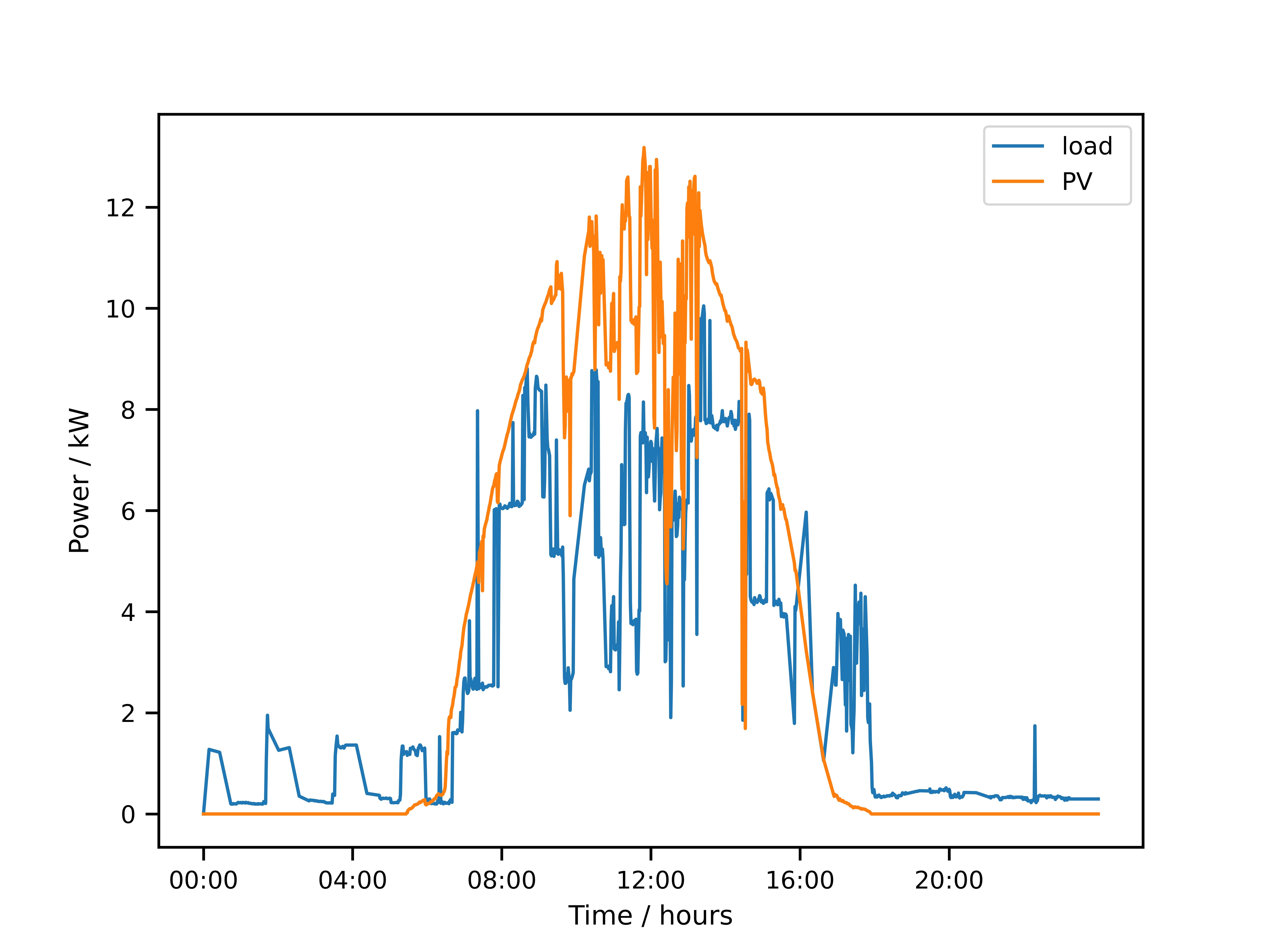}
    \caption{The load and PV production for the 1st April of 2023}
    \label{fig:rectangles2}
\end{figure}

Then, when we have our forecasting of when our device $i$ should start: ($\mathbf{T}_{device, ON_i}$), and when our device $i$ should end: ($\mathbf{T}_{device, OFF_i}$). We define an algorithm to be in concordance with the current production. 
This algorithm will thus adapt our original rectangles by creating different rectangles of different lengths (but with the same height) to adapt to the current production if the forecasted production is different (see Figure \ref{fig:rectangles2}).

\begin{algorithm}
\caption{Adaptive Device Algorithm}
\label{alg:adaptive_device}
\SetKwInput{Input}{Input}

\Input{
    $\mathbf{S}_{\text{device}}$: Device state (ON or OFF) \\
    $\mathbf{T}_{\text{device, ON}}$: Forecasted device ON time \\
    $\mathbf{T}_{\text{device, OFF}}$: Forecasted device OFF time \\
    $\mathbf{t}$: Current second of the day \\
    $\mathbf{N}$: Number of devices \\
    $\mathbf{T}_{\text{device, remaining}}$: Remaining ON time for each device
}

\For{$i \gets 1$ to $N$}{
    \eIf{$S_{\text{device}_i}$ is $OFF$}{
        \If{$t \geq T_{\text{device, ON}_i}$ and PV production is sufficient}{
            $S_{\text{device}_i} \gets ON$ \\
            $T_{\text{device, OFF}_i} \gets t + T_{\text{device, remaining}_i}$
        }
    }{
        \If{$t \geq T_{\text{device, OFF}_i}$}{
            $S_{\text{device}_i} \gets OFF$ \\
            $T_{\text{device, remaining}_i} \gets 0$
        }
        \If{$t < T_{\text{device, OFF}_i}$ and PV production is insufficient}{
            $S_{\text{device}_i} \gets OFF$ \\
            $T_{\text{device, remaining}_i} \gets T_{\text{device, OFF}_i} - t$
        }
    }
}
\end{algorithm}

\newpage

\section{Simulations}
\subsection{Time decision}
\label{sec:1}
We wanted to see the influence of the time frame for making decisions.
For that, we implemented a device management algorithm that works with decisions based on seconds, 5 minutes, 15 minutes, and 1 hour.

The algorithm works with the mode ‘adaptative’.
In the adaptive mode, we implement a dynamic threshold approach to optimise energy usage based on the PV production. A predefined threshold value is set, such as 4 kWs, which represents the maximum power consumption allowed from the PV production.
To illustrate the adaptive mode, let's consider a scenario where the PV production is 3.9 kWs at 1 pm, 4.1 kWs at 2 pm, 4.5 kWs at 3 pm, and 3.8 kWs at 4 pm. At each time interval, we compare the PV production with the threshold value.
For instance, at 1 pm, the PV production of 3.9 kWs is below the threshold of 4 kWs, so no device will be activated during this period. However, at 2 pm, when the PV production reaches 4.1 kWs, which exceeds the threshold, we initiate the activation of a controllable device.
Once the device is turned ON at 2 pm, it will consume the power from the PV production, in this case, 4 kWs per second, for the duration of the activation period. In our example, the device will remain ON until 4 pm, utilising the available 4 kWs per second during this time frame.
Finally, at 4 pm, when the PV production drops below the threshold to 3.8 kWs, we deactivate the device, turning it OFF as it can no longer operate at the desired power consumption level.

This adaptive mode allows us to dynamically adjust device activation based on real-time PV production, ensuring optimal energy utilisation within the constraints of the available power. By aligning device activation with higher PV production periods, we maximise the utilisation of clean energy and minimise reliance on the grid or other energy sources.
The adaptative mode can make devices work in parallel and can combine multiple devices. For example, for every second, we can say that the car takes 3 kWs, and the pool pump consumes 2 kWs. 

To measure the difference between 15-minute, 15-minute, and hour-based decisions, we decided to launch an adaptive mode with a fixed threshold of 4 kWs. 

\subsection{Comparison of device management algorithms}
To assess the efficacy of our device management algorithm, we conducted three simulations spanning from January 8th (considering the need for seven previous days) to April 10th, 2023, using 9 kWp solar panels.

The first simulation employed the smart device management algorithm described in this report, based on both current and forecasted PV production. Our automatic system prioritise the three controllable devices in the following order: the pool pump, the domestic hot water, and the car. Device A is considered to have higher priority than device B if, under the same conditions for activation, device A is switched ON first. Similarly, if given the same conditions for deactivation, device A is switched OFF last.
\label{sec:bruteforce}

In the second simulation, we scheduled the activation of the pool pump at 12 pm, domestic hot water at 7 pm, and the electric car at 6 pm each day. Each device was switched off after running for the same duration as in the first simulation. We chose these specific starting times based on our analysis of load data before the installation of the solar panels. This device management strategy is referred to as "bruteforce". This automatic system is independent of PV production.

\label{sec:adaptive}
The third simulation utilised the adaptive mode outlined in section \ref{sec:1}. In this mode, the threshold for each device's activation is based on its consumption level. This automatic system is based on the current PV production. 

It should be noted that the car was installed on March 25th, and, therefore, its consumption did not need to be managed before that date.

%% file: chapters/results.tex
\chapter{Results}

\section{Time decision}
\begin{figure}[H]
    \centering
    \includegraphics[width=11.5cm]{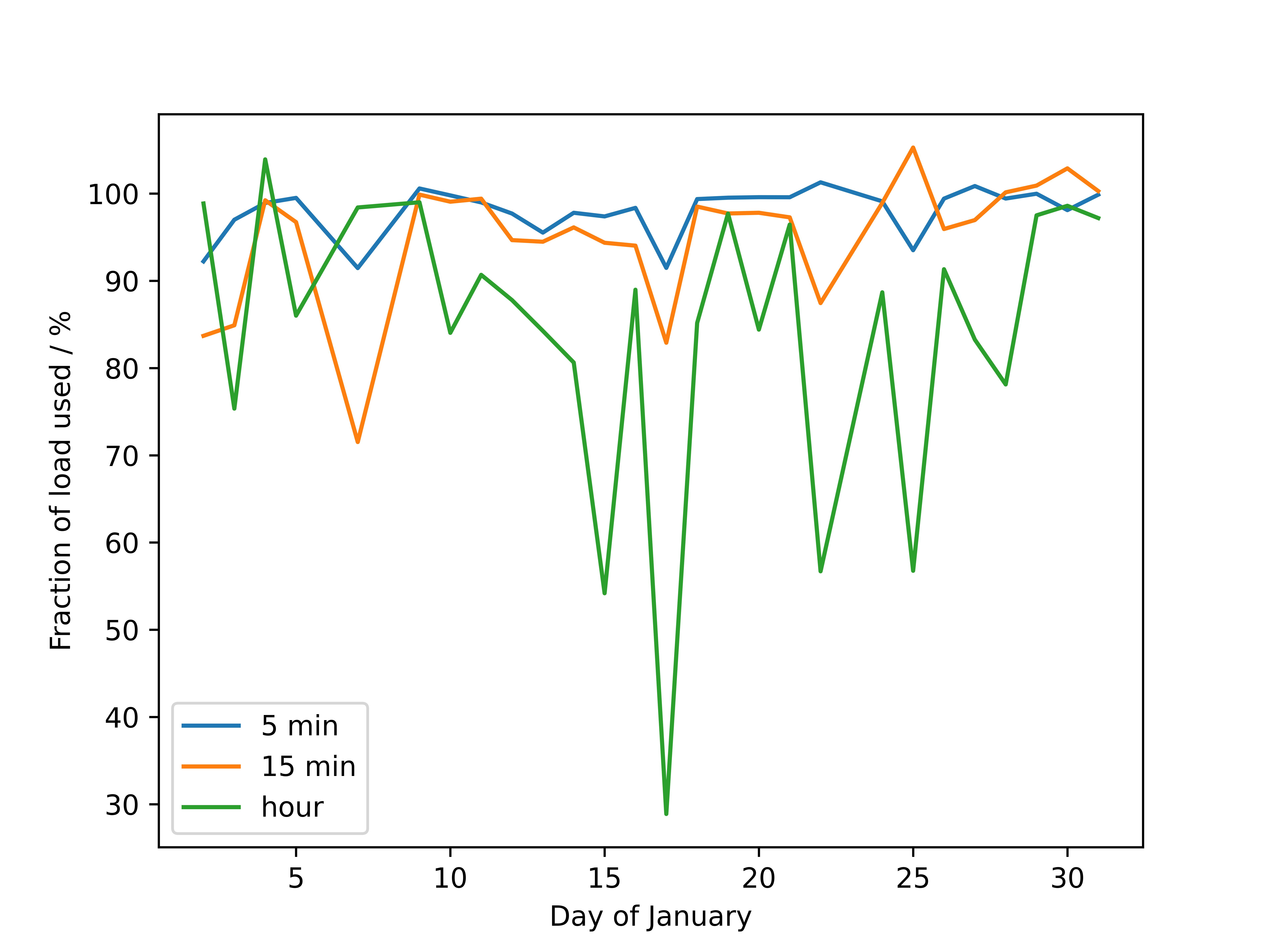}
    \caption{Proportion of energy used from PV production relative to the real-time decision using adaptive mode at a given 4 kW threshold for each day of January 2023}
    \label{fig:energyUsed}
\end{figure}

In order to evaluate the disparity in energy utilisation between decisions made at different time intervals (5 minutes, 15 minutes, and 1 hour) using the adaptive algorithm, we present the percentage of energy consumption relative to real-time decisions (measured in seconds) for each day throughout January 2023. Due to network connectivity issues and outages, we didn't have enough data to make the measurements desired with a high level of confidence on the 1st, 6th, 8th, and 23rd of January.
The results are shown in Figure \ref{fig:energyUsed}.
\begin{table}[H]
\centering
\begin{tabular}{||c|c|c|c||} 
 \hline
  & 5min & 15min & 1hour \\ [0.5ex] 
 \hline\hline
 mean & 98\% & 95\% & 84\% \\ 
 \hline
 median & 99\% & 97\% & 88\% \\
 \hline
 min & 91\% & 71\% & 28\% \\
 \hline
\end{tabular}
\caption{Statistics on the proportion of energy used from PV production relative to the real-time decision using adaptive mode at a given 4 kW threshold for each day of January 2023}
\label{tab:energyUsed}
\end{table}

\begin{figure}[H]
    \centering
    \includegraphics[width=17cm]{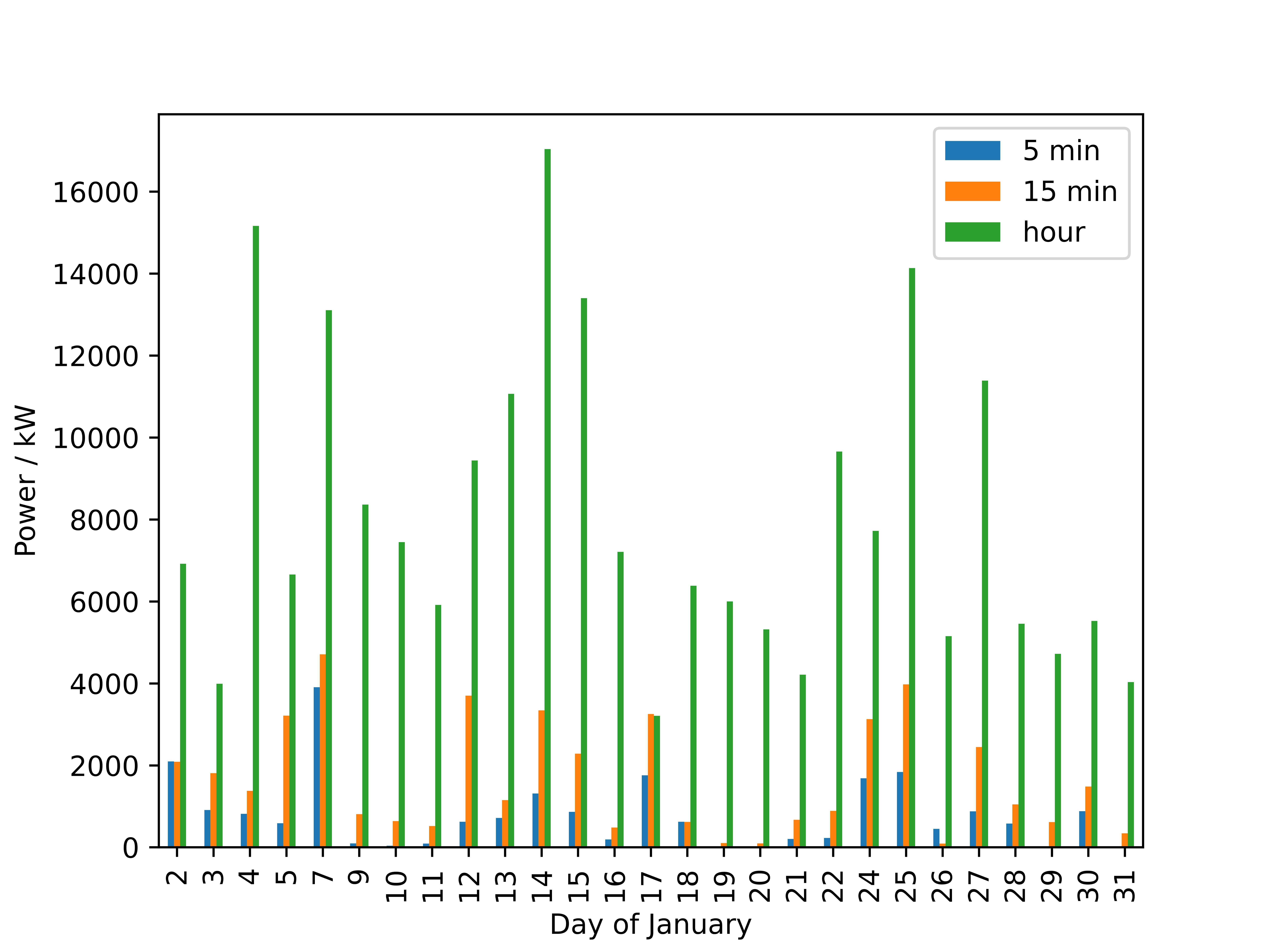}
    \caption{Energy paid from grid relative to the real-time decision using adaptive mode at a given 4 kW threshold for each day of January 2023}
    \label{fig:energyPaid}
\end{figure}

We notice that on the worst day of January 2023, an hour-based decision is using \textbf{28\%} of the available production from PV (\ref{tab:energyUsed}).
Then, we wanted to see how much we will pay for energy from the grid to satisfy the threshold of 4 kW, at different time intervals (5 minutes, 15 minutes, and 1 hour), relative to real-time decisions (measured in seconds) for each day throughout January 2023 (\ref{fig:energyPaid}).

\begin{table}[H]
\centering
\begin{tabular}{||c|c|c|c||} 
 \hline
  & 5min & 15min & 1hour \\ [0.5ex] 
 \hline\hline
 mean & 714.237kW & 1495.218kW & 7286.616kW \\ 
 \hline
 median & 581.408kW & 966.566kW & 6518.619kW \\
 \hline
 max & 3905.480kW & 4708.951kW & 17035.819kW \\
 \hline
\end{tabular}
\caption{Statistics on the energy paid from grid relative to the real-time decision using adaptive mode at a given 4 kW threshold for each day of January 2023}
\label{tab:energyPaid}
\end{table}

We notice that on the worst day of January 2023, hour-based decisions will cause us to pay 17035.819 kWs from the grid (\ref{tab:energyPaid}). 1 kWh costs 0.2€ in France \cite{energyprice}. 
17035.819 kWs equals approximately 4.73 kWh which makes approximately \textbf{0.95€} for the worst day of January 2023.

Over the first month period of real data, taking second-level decisions outperformed hourly decisions during 29 days. Over the whole period, taking second-based decisions saved \textbf{16\%} of energy compared to hour-based decisions.
In the following simulations of the report, we base our study purely on decisions by second.

\section{Comparison of device management algorithms}
\subsection{Comparison of smart and bruteforce management simulation}
\begin{figure}[H]
    \centering
    \includegraphics[width=10cm]{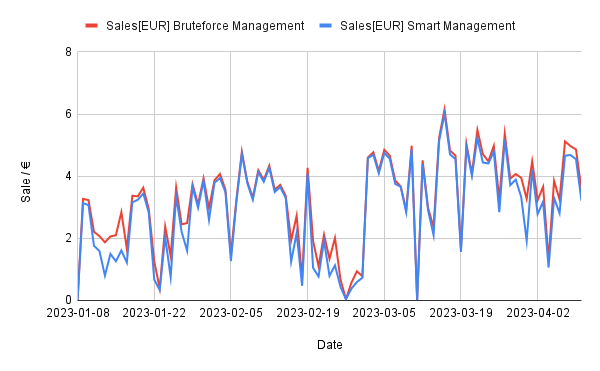}
    \caption{Sales of smart management and bruteforce management algorithms}
    \label{fig:sale_bruteforce}
\end{figure}

The simulation results presented in Figure \ref{fig:sale_bruteforce} depict that the bruteforce device management algorithm consistently achieves higher daily sales compared to the smart management algorithm. This disparity can be attributed to the fact that the smart management algorithm prioritises maximising the consumption of solar production overselling it, considering that the cost of purchasing energy exceeds the price of selling energy.

\begin{figure}[H]
    \centering
    \includegraphics[width=10cm]{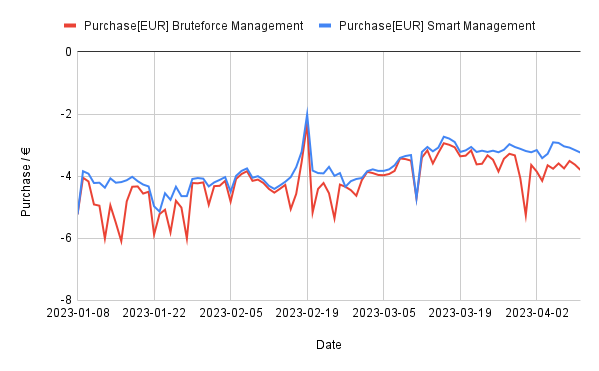}
    \caption{Purchases of smart management and bruteforce management algorithms}
    \label{fig:purchase_bruteforce}
\end{figure}

The concept of maximising the utilisation of PV production is further substantiated by Figure \ref{fig:purchase_bruteforce}, which demonstrates that the bruteforce device management algorithm incurs greater daily purchases compared to the smart management algorithm.

\begin{figure}[H]
    \centering
    \includegraphics[width=10cm]{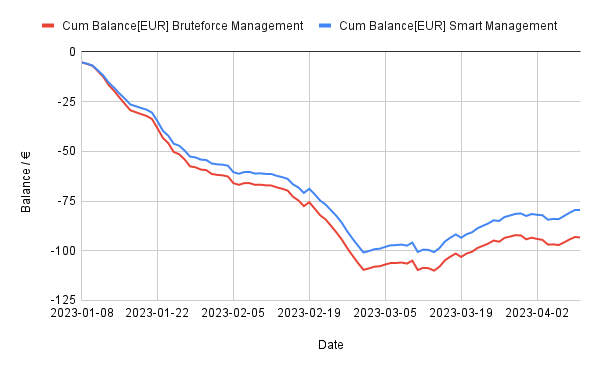}
    \caption{Cumulative balance of smart management and brute-force management algorithms}
    \label{fig:balance_bruteforce}
\end{figure}

To assess the overall performance of each device management approach, we compute the daily balance by summing the sales and purchases. Figure \ref{fig:balance_bruteforce} displays the cumulative balance for each management algorithm throughout the given period. Notably, the cumulative balance for the smart management algorithm consistently surpasses that of the bruteforce management algorithm.

\begin{figure}[H]
    \centering
    \includegraphics[width=10cm]{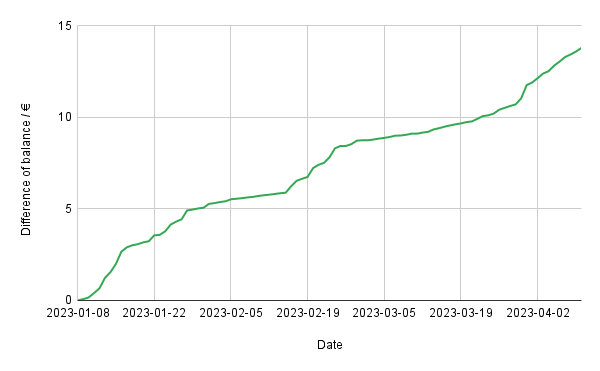}
    \caption{Difference of cumulative balance between smart management and bruteforce management algorithms}
    \label{fig:diffbalance_bruteforce}
\end{figure}

Furthermore, by analysing the statistics for the daily balance of each device management approach, we find that, on average, our algorithm achieves savings of approximately \textbf{0.43€} per day, which translates to approximately \textbf{157€} per year. Consequently, our algorithm yields \textbf{17.3\%} higher cost savings compared to the bruteforce approach.
However, it is important to note that the difference between the cumulative balance of the bruteforce management and the smart management varies over time, as depicted in Figure \ref{fig:diffbalance_bruteforce}, indicating a non-linear relationship.

In summary, our simulation-based comparison demonstrates the superior performance of the smart management algorithm in optimising energy consumption and financial savings. By maximising the utilisation of PV production, our algorithm consistently outperforms the bruteforce approach, resulting in significant cost reductions.

\subsection{Comparison of smart and Adaptative management simulation}
\begin{figure}[H]
    \centering
    \includegraphics[width=10cm]{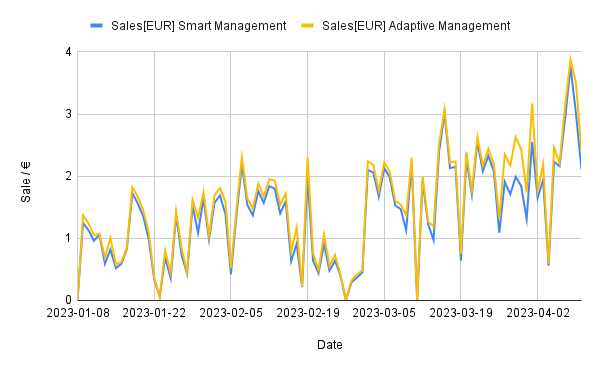}
    \caption{Sales of smart management and adaptive management algorithms}
    \label{fig:sale_adaptive}
\end{figure}

\begin{figure}[H]
    \centering
    \includegraphics[width=10cm]{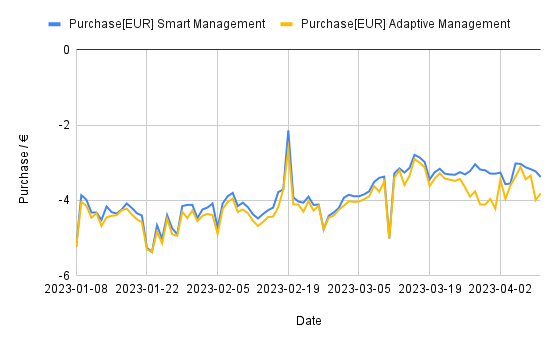}
    \caption{Purchases of smart management and adaptive management algorithms}
    \label{fig:purchase_adaptive}
\end{figure}

Next, we proceed to compare the results of the first simulation (Smart management) with the third simulation (Adaptive management). Similar to the previous comparison, the outcomes presented in Figure \ref{fig:sale_adaptive} demonstrate that the adaptive device management algorithm consistently achieves higher daily sales compared to the smart management algorithm. This can be attributed to the fact that the smart management algorithm possesses an optimal solution for maximising PV consumption compared to the adaptive management algorithm.
The notion of maximising the utilisation of PV production is further validated by Figure \ref{fig:purchase_adaptive}, which illustrates that the daily purchases made by the adaptive device management algorithm surpass those of the smart management algorithm.

\begin{figure}[H]
    \centering
    \includegraphics[width=10cm]{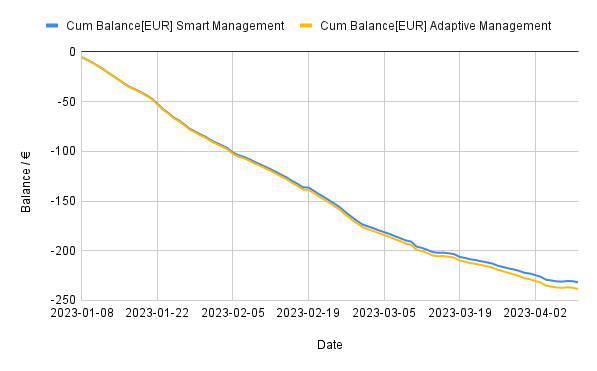}
    \caption{Cumulative balance of smart management and adaptive management algorithms}
    \label{fig:balance_adaptive}
\end{figure}

As in the previous comparison, we calculate the daily balance by summing the sales and purchases, and the results are shown in Figure \ref{fig:balance_adaptive}. Notably, the cumulative balance for the smart management algorithm consistently surpasses that of the adaptive management algorithm.

\begin{figure}[H]
    \centering
    \includegraphics[width=10cm]{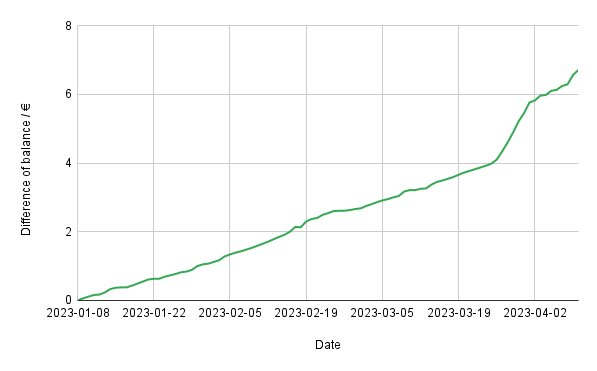}
    \caption{Difference of cumulative balance between smart management and adaptive management algorithms}
    \label{fig:diffbalance_adaptive}
\end{figure}

Moreover, by analysing the statistics for the daily balance of each device management approach, we observe an average savings of approximately \textbf{0.07€} per day, equivalent to approximately \textbf{26€} per year. Consequently, our algorithm achieves cost savings that are \textbf{2.8\%} higher than those of the adaptive management algorithm.
The difference in cumulative balances between the two device management algorithms is illustrated in Figure \ref{fig:diffbalance_adaptive}.

In summary, our comparative simulation of smart and adaptive device management algorithms demonstrates that the smart management algorithm outperforms the adaptive management algorithm in maximising energy consumption and achieving financial savings. By optimising PV utilisation, our algorithm consistently yields higher sales and cumulative balance, resulting in significant cost reductions compared to the adaptive management approach.

%% file: conclusion/conclusion.tex
\chapter{Conclusion}
\section{Summarisation}
In conclusion, this project has made significant advancements in energy forecasting and management. Firstly, in Chapter 2, we established the consumption time of each device, laying the foundation for further analysis.

We then succeed to forecast the solar production without any exogenous data by using a statistical method on the data of the seven past days. Our forecast of the PV production overestimates the real PV production per second with a mean of approximately \textbf{2.12\%} from the 8th of January to the 10th of April 2023.

In Chapter 4, we succeeded to find the consumption of each controllable device (the electric car, the pool pump, and the domestic hot water) for every day from the 8th of January to the 10th of April 2023 by using the median of the results of the seven past days for each day of our algorithm defined in section \ref{sec:algo_load} without using the "calm periods" defined in section \ref{sec:calm_periods}. 
We found a way to define the consumption of the other devices by defining the "base load" in section \ref{sec:base_load}.

In Chapter 5, we proved that real-time energy measures results in \textbf{2\%} higher daily charge (in kWh) of the car compared to 5 minutes resolution, \textbf{16\%} improvement for hourly measurement and reaction for the month of January 2023.
We furthermore defined optimal fit algorithms to shift the loads into the solar production envelope we forecasted in Chapter 3. We didn't forget to adapt our forecasting loads with the current real data per second.

Comparing the brute force management (section \ref{sec:bruteforce}) with the smart management algorithm developed in this project, we found that the smart management approach resulted in \textbf{17.3\%} more cost savings from the 8th of January to the 10th of April 2023. This result can be compared to the study conducted by Paetz et al., where they observed a cost saving of 17\% during their experiment involving a semi-automatic shift-loading system in a household located in Germany \cite{paetz2013load}. Furthermore, when comparing the adaptive management (section \ref{sec:adaptive}) with the smart management algorithm, the latter achieved an \textbf{2.8\%} in cost savings over the same period.

These results highlight the effectiveness of our forecasting and management approaches in optimising energy consumption and cost savings. The findings serve as a valuable contribution to the field of energy management and lay the groundwork for future research and implementation in real-world scenarios.

\section{Limitations}
This report is subject to certain limitations that should be acknowledged. Firstly, the data used for analysis and evaluation were limited to the period from the 1st of January to the 10th of April 2023. As a result, the findings and conclusions may not fully capture the variations and patterns that could emerge during the summer and autumn seasons. Future studies should consider expanding the data collection period to include a more comprehensive representation of seasonal changes.

Secondly, the analysis conducted in this report was based solely on the data collected from a single household located in the South of France. Although this dataset provided valuable insights into the energy consumption patterns and solar power generation of this specific household with a 14.7 kWc PV installation, it may not be fully representative of households in other regions or with different PV system configurations. Therefore, caution should be exercised when generalising the findings to broader contexts.

To address these limitations and ensure the robustness and applicability of future research, it is recommended to gather data from a more diverse range of households in different geographical locations, with varying PV installation capacities. By incorporating a broader dataset, researchers can gain a more comprehensive understanding of the factors influencing energy consumption and solar power generation across different contexts.

Moreover, extending the data collection period to include all four seasons would provide a more comprehensive analysis of how energy consumption and solar power generation patterns evolve throughout the year. This expanded timeframe would allow for a better understanding of the seasonal variations and their impact on the proposed forecasting algorithms and energy management strategies.

\section{Extensions}
To extend the work accomplished in this project, our objective is to develop an algorithm capable of accurately forecasting the energy consumption patterns of controllable devices, without relying on user input. This holds significant importance as it can lead to optimised energy usage, cost reduction, and improved energy efficiency in households worldwide.

One of the key challenges in achieving this goal is to create an algorithm that can adapt to any household, irrespective of the specific devices used. By eliminating the need for defining the consumption time of each controllable device, we can ensure the algorithm's applicability in diverse settings.

Furthermore, it is essential to enhance the accuracy of the energy production forecast by considering external factors such as weather conditions. While our study utilised data from a household in the South of France, where many days were sunny and closely aligned with the clear sky model, incorporating numerical weather predictions (NWPs) can significantly improve the forecast accuracy.

NWPs, including factors such as Direct Normal Irradiance (DNI), Diffuse Horizontal Irradiance (DHI), and Global Horizontal Irradiance (GHI), play a vital role in solar power forecasting. By integrating NWPs as inputs into the forecasting algorithm, we can account for variations in solar radiation and adjust the forecasted production accordingly.
The paper \cite{bacher2009online} used these NWPs as input and achieved an improvement of around 35\% with his ARX model. DNI, DHI and GHI can be estimated using the REST2 model described in this paper \cite{gueymard2008rest2}.